\documentclass[12pt]{article}
\usepackage{ifthen}

\newboolean{tab}
\setboolean{tab}{false}

\usepackage[latin1]{inputenc}
 \usepackage[authoryear]{natbib}
\usepackage{graphicx}
\usepackage{latexsym}
\usepackage{amssymb}
\usepackage{epsfig}
\makeatletter\@addtoreset{equation}{section}\makeatother

\def\eps{\varepsilon}
\def \R{\mathbb{R}}
\def \E{\mathbb{E}}

\def \Z{\mathbb{Z}}

\newcommand{\Pe}{\mathbb{P}}

\newcommand{\prodi}{\mathop{\mbox{\Huge{$\pi$}}}}

\newcommand{\Dkonv}{\stackrel{\mathcal{D}}{\rightarrow}}
\newcommand{\Pkonv}{\stackrel{P}{\rightarrow}}

\newcommand{\G}{\mathcal{G}}
\newcommand{\D}{\mathcal{D}}

\newcommand{\M}{\mathcal{M}}

\newcommand{\Ne}{\mathcal{N}}

\newcommand{\T}{\mathcal{T}}
\newcommand{\B}{\mathcal{B}}
\newcommand{\Rr}{\mathcal{R}}
\newcommand{\Zz}{\mathcal{Z}}

\newcommand{\Pp}{\mathcal{P}}

\newcommand{\bb}{\mathbf{b}}
\newcommand{\hh}{\mathbf{h}}

\newcommand{\ZZ}{\mathbf{Z}}

\newcommand{\WWW}{\mathbb{W}}
\newcommand{\VVV}{\mathbb{V}}

\newcommand{\bean}{\begin{eqnarray*}}
\newcommand{\eean}{\end{eqnarray*}}
\newcommand{\bea}{\begin{eqnarray}}
\newcommand{\eea}{\end{eqnarray}}
\newcommand{\be}{\begin{eqnarray}}
\newcommand{\ee}{\end{eqnarray}}
\newcommand{\beq}{\begin{equation}}
\newcommand{\eeq}{\end{equation}}

\newtheorem{theo}{Theorem}[section]
\newtheorem{lemma}[theo]{Lemma}
\newtheorem{cor}[theo]{Corollary}
\newtheorem{rem}[theo]{Remark}

\begin{document}

\title{Censored quantile regression processes under dependence and penalization.}

\author{Stanislav Volgushev$^{\normalsize\rm a,b}$\thanks{
Supported in part by the Collaborative Research Center ``Statistical modeling of nonlinear dynamic processes'' (SFB 823, Teilprojekt C1) of the German Research Foundation (DFG) and by the DFG grant Vo1799/1-1. Part of this research was conducted while the first author was visiting the University of Illinois at Urbana-Champaign. He would like to thank the people at the Statistics and Economics departments for their hospitality.
The authors would like to thank Roger Koenker and Steve Portnoy for many helpful hints and discussions. Our thanks also go to Judy Wang for sending us a preprint of her paper.}, Jens Wagener, Holger Dette$^{\normalsize\rm a*}$
\\
$^{\normalsize\rm a}$ Ruhr-Universit\" at Bochum \vspace{1mm}
\\
$^{\normalsize\rm b}$ University of Illinois at Urbana-Champaign.
$\;$\vspace{-7mm}\\
}
\date{}
\maketitle

\begin{abstract}
We consider quantile regression processes from censored data under dependent data structures and derive a uniform Bahadur representation for those processes. We also consider cases where the dimension of the parameter in the quantile regression model is large. It is demonstrated that traditional penalization methods such as the adaptive lasso yield sub-optimal rates if the coefficients of the quantile regression cross zero. New penalization techniques are introduced which are able to deal with specific problems of censored data and yield estimates with an optimal rate. In contrast to most of the literature, the asymptotic analysis does not require the assumption of independent observations, but is based on
rather weak assumptions, which are satisfied for many kinds of dependent data.
\end{abstract} 



\section{Introduction} \label{intro}

Quantile regression for censored data has found considerable attention in the recent literature. Early work dates back to
\cite{powell1984}, \cite{powell1986} and \cite{newepowe1990} who proposed quantile regression methods in the case where all censoring variables are known [see also \cite{fitzenberger1997}]. \cite{yingetal95} introduced median regression in the presence of right independent censoring. Similar ideas were considered by \cite{BangTsiatis02} and later \cite{zhou06}.

All these papers have in common that the statistical analysis requires the independence of the censoring times and covariates.
\cite{portnoy2003} and \cite{portlin2010} replaced this rather strong assumption by conditional independence of survival and censoring times conditional on the covariates. The resulting iterative estimation procedure was based on the principle of mass redistribution that dates back to the Kaplan-Meier estimate. An alternative and very interesting quantile regression method for survival data subject to conditionally independent censoring was proposed by \cite{penghuang08} and \cite{huang2010} who exploited an underlying martingale structure of the data generating mechanism. In particular, in the four last-mentioned papers weak convergence of quantile processes was considered. This is an important question since it allows to simultaneously analyze the impact of covariates on different regions of the conditional distribution. We also refer to the recent work of \cite{wangwang2009}, \cite{lengtong2012} and \cite{tanwanhezhu2012} who discussed quantile regression estimates that cope with censoring by considering locally weighted distribution function estimators and employing mass-redistribution ideas. All of the references cited above have in common that the asymptotic analysis is rather involved and relies heavily on the assumption of independent observations. An important and natural question is, whether, and how far, this assumption can be relaxed. One major purpose  of the present paper   is to demonstrate that a sensible asymptotic theory can be obtained under rather weak assumptions on certain empirical processes that are satisfied for many kinds of dependent data.
We do so by deriving a uniform Bahadur representation for the quantile process. In some cases, we also discuss the rate of the remainder term.

The second objective of this paper deals with settings where the dimension of the parameter of the quantile regression model is large. In this case the estimation problem  is intrinsically harder. Under sparsity assumptions  penalized   estimators can yield substantial improvements in estimation accuracy. At the same time, penalization allows to identify those components of the predictor which have an impact on the response. In the uncensored case, penalized quantile regression has found considerable interest in the recent literature [see \cite{zouyuan2008}, \cite{wuliu2009} and \cite{belche2011} among others].
On the other hand - to the best knowledge of the authors - there are only three papers which discuss penalized estimators  in the context of censored quantile regression. \cite{showzhan2010} proposed to penalize the estimator developed in \cite{zhou06} by an adaptive lasso penalty. These authors assumed unconditional independence between survival and censoring times and considered only the median.
\cite{wangzhouli2012} proposed to combine weights that are estimated by local smoothing with an adaptive lasso penalty. The authors considered a model selection at a fixed quantile and did not investigate process convergence of the corresponding estimator.

In contrast to that, \cite{wagvoldet2012} investigated sparse quantile regression models and properties of the quantile process in the context of censored data. As \cite{showzhan2010}, these authors assumed independence of the censoring times and predictors, which may not be a reasonable assumption in many practical problems and moreover might lead to inefficient estimators [see the discussion in \cite{koenker2008} and \cite{portnoy2009}]. An even more important point reflecting the difference between the philosophy of quantile versus mean regression   was not considered in the last-named paper. In contrast to   mean, quantile regression is concerned with the impact of predictors on different parts of the distribution. This implies that the set of important components of the predictor could vary for different quantiles. For example, it might be possible that a certain part of the predictor has a strong influence on the $95\%$-quantile of the distribution of the response, while a different set relates to the median. Also, quantile coefficients might cross zero as the probability for which the quantile regression is estimated varies. Traditional analysis of penalized estimators, including the one given in \cite{wagvoldet2012}, fails in such situations. At the same time, it might not be reasonable to exclude covariates from the model just because they have zero influence at a fixed given quantile. All those considerations demonstrate the need for penalization techniques that take into account the special features of quantile regression. To the best of our knowledge,  no results answering these questions are available in the context of censored quantile regression.

Therefore the second purpose of the present paper is to construct novel penalization techniques that are flexible enough to deal with the particular properties of censored quantile regression, and to provide a rigorous analysis of the resulting quantile regression processes. One major challenge for the theoretical analysis of censored regression quantiles in the present setting is the sequential nature of the underlying estimation procedures. While in other settings  estimators for different quantiles do not interact, the situation is fundamentally different in the case of censored data when iterative procedures need to be applied. In the course of our analysis, we demonstrate that using
traditional  generalizations of concepts from the mean regression setting can result in sub-optimal rates of convergence. As a solution of this problem we propose penalties that avoid this problem and additionally allow to analyze the impact of predictors on quantile regions instead of individual quantiles. Finally, all our results hold for a wide range of dependence structures thus considerably extending the scope of their applicability.\\
\\
The remaining part of the paper is organized as follows. The basic setup is introduced in Section \ref{assu}. In Section \ref{dependent}, we concentrate on the properties of the unpenalized estimator in settings where the realizations need not be independent and derive a uniform Bahadur representation. Various ways of penalizing the censored quantile process and the properties of the resulting estimators are discussed in Section \ref{genpen}. A small simulation study illustrating the findings in this section is presented in Section \ref{simul}.
Finally, all proofs and technical details are deferred to an appendix in Section \ref{proofs}.


\section{Censored quantile regression} \label{assu}
We consider a censored regression problem with response $T_i$, predictor $\ZZ_i$ and censoring time $C_i$,
where the random variables $T_i$ and  $C_i$ may be dependent, but conditionally on the $d$-dimensional covariate
$\ZZ_i$ the response $T_i$ and the censoring time $C_i$ are independent. As usual we assume that
instead of $T_i$ we only observe $X_i = T_i\wedge C_i,$ and the indicator $\delta_i = I\{X_i = T_i\}.$  Let
$\{ T_i,C_i,\ZZ_i \}_{i=1}^n$ denote $n$ identically distributed copies of the random variable
$(T_1,C_1,\ZZ_1)$. The aim consists in statistical inference regarding the quantile function of the random variable $T$
conditional on the covariate vector $\ZZ$ on the basis of the sample $\{X_i, \ZZ_i, \delta_i\}^n_{i=1}$. In particular we would like
to study the influence of the components of the predictor on different quantiles of the distribution of $T$.
Following \cite{portnoy2003} and \cite{penghuang08}, we assume that the conditional quantile functions of $T$ are linear in $\ZZ$, i.e.
\beq \label{model}
Q_\tau(T|\ZZ) : = \inf\{t: P(T\leq t|\ZZ) \geq \tau\}   = \ZZ^t\beta(\tau)
\eeq
for $\tau \in [\tau_L,\tau_U] \subset [0,1)$. Combining ideas from the above references, an estimator for the coefficient function $\beta(\tau)_{\tau \in [\tau_L,\tau_U]}$ can be constructed in an iterative manner. To be precise, consider a uniformly spaced grid
\be \label{grid}
\tau_L <\tau_1<...<\tau_{N_\tau(n)}=\tau_U
\ee
with width $a_n=o(n^{-1/2})$ and set $b_n := a_n/(1-\tau_U)$. The estimator for $\beta(\tau)$ is now defined as a piecewise constant function.
We follow \cite{portnoy2003} by assuming that there is no censoring below the $\tau_L$'th quantile where $\tau_L > 0$. Setting $\tau_0 = \tau_L$, the estimator $\hat\beta(\tau_L)$ is defined as the classical \cite{koebas1978} regression quantile estimator without taking censoring into account. For $j=1,\dots,N_{\tau(n)}$ the estimator $\hat\beta(\tau_j)$ of $\beta(\tau_j)$ is then sequentially defined as any value from the set of minimizers of the   convex function
\beq \label{intesteq.ava_np}
\tilde H_j(\bb) := \frac{1}{n} \sum_i \Big( \delta_i|X_i - \ZZ_i^t\bb| - \ZZ_i^t\bb \Big( \delta_i -  2\int_{[\tau_0,\tau_j)}I\{X_i \geq \ZZ_i^t\hat\beta(u)\}dH(u) - 2\tau_0 \Big)\Big) \quad
\eeq
Here $H(u) := -\log(1-u)$ and $\hat \beta(\tau)$ is defined as constant and equal to $\hat \beta(\tau_j)$ whenever $\tau \in [\tau_j,\tau_{j+1})$.
The convexity of $H$ greatly facilitates the computation of the estimators. In particular the computation of the directional derivative of the function $\tilde H_j$ at the point $\bb$ in direction of $\xi$ yields
\bea \label{esteq.ava0_np}
\Psi_j(\bb,\xi) &=& \frac{- 2}{n}\sum_{i=1}^n \xi^t\ZZ_i \Big( N_i(\ZZ_i^t\bb) - \int_{[\tau_0,\tau_j)}I\{X_i \geq \ZZ_i^t\hat\beta(u)\}dH(u) - \tau_0 \Big)
\\ \nonumber
&&  + \frac{1}{n}\sum_{i=1}^n I\{X_i = \ZZ_i^t\bb\}(\delta_i\xi^t\ZZ_i + |\xi^t\ZZ_i|)
\eea
where $N_i(t) := \delta_i I\{X_i \leq t\}$ and $sgn(a) := \frac{a}{|a|}$ if $a\neq 0$ with $sgn(0) := 0$. We thus obtain that  any minimizer $\hat \bb$ of the function $H$ defined in (\ref{intesteq.ava}) satisfies the condition
\beq \label{esteq.ava_np}
\inf_\xi \Psi_j(\hat \bb,\xi) \geq 0.
\eeq
The first major contribution of the present paper consists in replacing the i.i.d. assumption that underlies all asymptotic investigations considered so far by general conditions on certain empirical processes. In particular, we demonstrate that these conditions are satisfied for a wide range of dependency structures. Moreover, instead of providing results on weak convergence, we derive a uniform (weak) Bahadur representation that can be used as starting point for the investigation of general L-type statistics [see e.g. \cite{porko1989}] and rank-based testing procedures [see \cite{gujukopo1993}].

\begin{rem} \label{rem:ph}
{\rm
\cite{penghuang08} studied a closely related estimate. More precisely these authors proposed to set $\hat\beta(0) := 0$ defined their estimator for $\beta(\tau_j)$ as the iterative (generalized) solution of the equations
\[
\sum_{i=1}^n \ZZ_i\Big(N_i(\ZZ_i^t\bb) - \int_0^{\tau_j} I\{X_i \geq \ZZ_i^t \hat\beta(u)\} dH(u)\Big) \approx 0
\]
Note that this corresponds to the first line in the definition of $\Psi_j$ in equation (\ref{esteq.ava0_np}). In the case when the $X_i$ have a continuous distribution, the second line in the definition of $\Psi_j$ is of order $O_P(1/n)$ uniformly with respect to $\bb$. Therefore (under this additional assumption) this part is negligible compared to the rest of the equation and the proposed estimator can thus be viewed as the solution of the estimating equation
\[
\sum_{i=1}^n \ZZ_i\Big(N_i(\ZZ_i^t\bb) - \int_{[\tau_0,\tau_j)} I\{X_i \geq \ZZ_i^t \hat\beta(u)\} dH(u) - \tau_0\Big) \approx 0
\]
which corresponds to the one considered by \cite{penghuang08} if we set $\tau_0 = 0$.
}
\end{rem}

\begin{rem} \label{rem:phport}
{\rm 
It is possible to show that in the case  with no censoring up to a quantile $\tau_L$, the estimator starting at $\tau_L$ and the version starting at $\tau_0=0$ considered by \cite{penghuang08} share the same limiting behavior. However, we would like to point out that, in order for the estimator starting at $\tau_0 = 0$ to be well-behaved, conditions controlling all the lower part of the conditional distribution of the survival time given the covariates need to be imposed. Obviously, no such assumptions are necessary for the version starting at $\tau_L$, and for this reason this version   seems to be preferable in cases where there is no censoring below a certain quantile.}
\end{rem}

It is well known that in models with insignificant coefficients penalization of the estimators can yield significant improvements in the estimation accuracy. At the same time, this method allows for the identification of the components of the predictor which correspond to the non-vanishing components of the parameter vector. The second part of our paper is therefore devoted to considering penalized versions of the estimator described above. Penalization can be implemented by adding an additional term to the estimating equation in (\ref{intesteq.ava_np}). More precisely, we propose to define
\[
\hat \beta(\tau_0) := \arg\min_\bb \frac{1}{n} \sum_i \rho_{\tau_0}(X_i - \bb^t\ZZ_i) + \lambda_n\sum_{k=1}^d |\bb_k|/p_k(n,\tau_0)
\]
and replace the function $\tilde H_j$ in (\ref{intesteq.ava_np}) by
\bea \label{intesteq.ava}
H_j(\bb) &:=& \frac{1}{n} \sum_i \Big( \delta_i|X_i - \ZZ_i^t\bb| - \ZZ_i^t\bb \Big( \delta_i -  2\int_{[\tau_0,\tau_j)}I\{X_i \geq \ZZ_i^t\hat\beta(u)\}dH(u) - 2\tau_0 \Big)\Big) \quad
\\ \nonumber
&& + 2\lambda_n\sum_{k=1}^d |\bb_k|/p_k(n,\tau_j).
\eea
Here, the quantity $p(n,\tau_j)=(p_1 (n, \tau_j),\dots,p_d(n,\tau_j))$ denotes a $d$-dimensional vector that, together with $\lambda_n$, controls the amount of penalization and is allowed to depend on the data. A very natural choice is given by a version of the adaptive lasso [see \cite{zou2006}], that is
\beq\label{adappen}
p_k(n,\tau_j) = |\tilde\beta_k(\tau_j)|
\eeq
($k=1,\ldots ,d$) where $\tilde\beta(\tau)$ is some preliminary estimator for the parameter $\beta(\tau)$. A detailed discussion of estimators based on this penalization will be given in Section \ref{detailsadapt}. In particular, we will demonstrate that in certain situations the adaptive lasso can lead to non-optimal convergence rates. Alternative ways of penalization that avoid this problem will be discussed in Section \ref{detailsint}.

\begin{rem}
{\rm Note that we also allow the choice $p_k(n,\tau_j) = \infty$ throughout this paper if it is not stated otherwise. By this choice   we do not use a penalization for the $k$'th component, which would be reasonable if a variable is known to be important. For example, it is reasonable not to penalize the component of $\beta$ corresponding to the intercept since it will typically vary across quantiles and thus be different from zero. }
\end{rem}

\section{A Bahadur representation for dependent data}\label{dependent}
For the asymptotic results, we will need the following notation and technical assumptions which are collected here for later reference.
Consider the conditional distribution functions
\begin{eqnarray*}
\tilde F(t|z) :&=& P(X\leq t|z),\quad F(t|z) := P(X \leq t, \delta=1|z)
\end{eqnarray*}
and denote by $\tilde f(t|z), f(t|z)$ the corresponding conditional densities. Define the quantities
\bea \label{mu}
\mu(\bb) &:=& \E[\ZZ I\{X\leq\ZZ^t\bb,\delta=1\}], \quad \tilde\mu(\bb) := \E[\ZZ I\{X \geq \ZZ^t\bb\}]
\\ \label{nu}
\nu_n(\bb) &:=& \frac{1}{n}\sum_i \ZZ_iN_i(\ZZ_i^t\bb) - \mu(\bb), \quad \tilde \nu_n(\bb) := \frac{1}{n}\sum_i \ZZ_i I\{X_i \geq \ZZ_i^t\bb\} - \tilde \mu(\bb).
\eea

We need the following conditions on the data-generating process.

\begin{itemize}
\item[(C1)] \label{c:z1} The model contains an intercept, that is $Z_{i,1} = 1$ a.s. for $i=1,...,n$ and there exists a finite constant $C_Z > 0$ such that $  \|\ZZ\| \leq C_Z$ a.s. [here and throughout the paper, denote by $\|\cdot\|$ the maximum norm].
\item[(C2)] \label{c:lbeta} There exist a finite constant $C_4$ such that $$\|\beta(\tau_1)-\beta(\tau_2)\| \leq C_4 |\tau_1-\tau_2|$$
for all $\tau_1,\tau_2 \in [\tau_L,\tau_U]$.
\item[(C3)] \label{c:Z} Define the set $\B(\T,\eps) := \{\bb \in \R^d: \inf_{\tau \in \T}\|\bb - \beta(\tau)\|<\eps\}$. Then
\[
\sup_{\bb \in \B(\T,\eps)} \sup_{z} f(z^t\bb|z) =: K_f < \infty, \quad \sup_{\bb \in \B(\T,\eps)} \sup_{z} \tilde f(z^t\bb|z) =: \tilde K_f < \infty
\]
Moreover $f, \tilde f$ are uniformly continuous on $\{\bb^tz: \bb \in \B(\T,\eps), z \in \Zz\}\times \Zz$ with respect to both arguments and uniformly H÷lder continuous with respect to the first argument, i.e. for some $\gamma>0$ and $H_f,\tilde H_f < \infty$
\bean
&&\sup_{\bb_1,\bb_2 \in \B(\T,\eps)} \sup_{z} |f(z^t\bb_1|z) - f(z^t\bb_2|z)| \leq H_f\|\bb_1 - \bb_2\|^\gamma,
\\
&& \sup_{\bb_1,\bb_2 \in \B(\T,\eps)} \sup_{z} |\tilde f(z^t\bb_1|z) - \tilde f(z^t\bb_2|z)| \leq \tilde H_f\|\bb_1 - \bb_2\|^\gamma
\eean
\item[(C4)] We have
\[
\inf_{\bb \in \B(\T,\eps)}\lambda_{\min}(\E[(\ZZ\ZZ^tf(\ZZ^t\bb|\ZZ)]) =: \lambda_0 >0
\]
where $\lambda_{\min}(A)$ denotes the smallest eigenvalue of the matrix $A$.
\end{itemize}

\medskip

\begin{rem}
{\rm
Condition (C1) has been imposed by all authors who considered model (\ref{model}). While it possibly could be relaxed, this would introduce additional technicalities and we therefore leave this question to future research. Conditions (C2),(C3) place mild restrictions on the regularity of the underlying data structure. Condition (C4) is similar to condition (C4) in \cite{penghuang08}. It yields an implicit characterization of the largest quantile that is identifiable in the given censoring model. For a more detailed discussion of this point, we refer the interested reader to Section 3 of \cite{penghuang08}. 
}
\end{rem}

In contrast to most of the literature in this context which requires independent observations, our approach is based on a general condition on certain empirical processes which holds for many types of dependent data. More precisely, we assume the following conditions.
\begin{enumerate}
\item[(D1)] With the notation (\ref{nu}) we have
\[
\sup_{\bb \in \R^d} \|\nu_n(\bb)\| + \sup_{\bb \in \R^d} \|\tilde\nu_n(\bb)\| = o_P(1) 
\]
\item[(D2)] For some $\eps >0$ define $\B := \{\bb: \inf_{\tau \in [\tau_L,\tau_U]}\|\bb- \beta(\tau)\| \leq \eps\}$ and for a function $g$ on $\B$ define
\[
\omega_{a}(g) := \sup_{\|\bb_1-\bb_2\|\leq a, \bb_1,\bb_2\in\B} \|g(\bb_1)-g(\bb_2)\|
\]
Then the empirical processes $(\sqrt{n}\nu_n(\bb))_{\bb\in\B}$ and $(\sqrt{n}\tilde\nu_n(\bb))_{\bb\in\B}$ satisfy for any $a_n=o(1)$
\[
\omega_{a_n}(\sqrt{n}\nu_n) = o_P(1), \quad \omega_{a_n}(\sqrt{n}\tilde\nu_n) = o_P(1).
\]
\item[(D3)] The process
\[
w_n(s) := \frac{\tau_0}{n}\sum_{i=1}^n (\ZZ_i - \E\ZZ_i)- \nu_n(\beta(s)) + \int_{[\tau_0,s)}\tilde \nu_n(\beta(u))dH(u)
\]
indexed by $s \in [\tau_L,\tau_U]$ converges weakly towards a centered Gaussian process $\WWW$.
\end{enumerate}

First of all, we would like to point out that for independent data, conditions (D1)-(D3) follow under (C1) and (C3) and in this case
\[
\omega_{a_n}(\sqrt{n}\nu_n) + \omega_{a_n}(\sqrt{n}\tilde\nu_n) = O_P((a_n\log n)^{1/2}\vee (n^{-1}\log n)^{1/2}).
\]
We now provide a detailed discussion of results available in settings where the independence assumption is violated. To this end, note that $\nu_{n,k}(\bb) = \int g_\bb dP_n - \E[g_\bb(\ZZ,X,\delta)]$ where $g_\bb(z,x,\delta) := z_k I\{x\leq z^t\bb\}\delta$ and $P_n$ denotes the empirical measure of the observations $(X_i,\ZZ_i,\delta_i)_{i=1,...,n}$. Thus for any set $B \subset \R^d$ the process $(\sqrt{n}\nu_{n,k}(\bb))_{\bb\in B}$ can be interpreted as empirical process indexed by the class of functions $\{g_\bb|\bb\in B\}$.

\begin{rem} \label{rem:vc}
{\rm
Combining Lemma 2.6.15 and Lemma 2.6.18 from \cite{vaarwell1996} shows that $\{g_\bb|\bb\in \R^d\}$ is VC-subgraph [see Chapter 2.6 in the latter reference for details], and under assumption (C1) all functions in this class are uniformly bounded. Similar arguments apply to $\tilde \nu_{n,k}(\bb)$. The problem of uniform laws of large numbers for VC-subgraph classes of functions for dependent observations has been considered by many authors. A good overview of recent results can be found in \cite{adanobe2010} and the references cited therein. In particular, the results in the latter reference imply that (D1) holds as soon as $(X_i,\ZZ_i,\delta_i)_{i \in \Z}$ is ergodic, (C1) is satisfied and the conditional distribution function of $X$ given $\ZZ$, i.e. $\tilde F$, is uniformly continuous with respect to the first argument.
}
\end{rem}

\begin{rem} \label{rem:wk}
{\rm
Condition (D2) essentially imposes uniform asymptotic equicontinuity of the processes $n^{1/2}\nu_n$, $n^{1/2}\tilde\nu_n$. It is intrinsically connected to weak convergence of those processes. More precisely, Theorem 1.5.7, Addendum 1.5.8 and Example 1.5.10 in \cite{vaarwell1996} imply that (D2) will hold as soon as the processes $n^{1/2}\nu_n$, $n^{1/2}\tilde\nu_n$ converge weakly towards centered Gaussian processes, say $\VVV, \tilde \VVV$, with the additional property that $\E[(\VVV(\bb_1)-\VVV(\bb_2))^2] = o(1)$ implies $\|\bb_1-\bb_2\| = o(1)$. Condition (D2) can thus be checked by establishing weak convergence of $n^{1/2}\nu_n$, $n^{1/2}\tilde\nu_n$ and considering the properties of their covariance. The literature on weak convergence of processes indexed by certain classes of functions in dependent cases is rather rich.\\
Specifically, with the notation from Remark \ref{rem:vc}, it is possible to show that under assumption (C3) the bracketing numbers [see Definition 2.1.6 in \cite{vaarwell1996}] of the class $\G := \{g_\bb|\bb\in \B(\T,\eps)\}$ satisfy $\Ne_{[~]}(\eps,\G,P_{X,\ZZ,\delta}) \leq c\eps^{-d}$ for some finite constant $c$. Thus, among many others, the results from \cite{aryu1994} for $\beta$-mixing, the results from \cite{anpo1994} for $\alpha$-mixing and the results from \cite{hagemann2012stochastic} for data from general non-linear time series models can be applied to check condition (D2). For example, the results in \cite{aryu1994} imply that (D2) will hold as soon as $(\ZZ_i,T_i,R_i)_{i\in\Z}$ is a strictly stationary, $\beta-$mixing sequence with coefficients $\beta_k = O(k^{-r})$ for some $r>1$.
}
\end{rem}

We now are ready to state the main result of this section.

\begin{theo} \label{bahadur}
Assume that $\tau_0 = \tau_L >0$, that for some $a>0$ we have $P(C > \ZZ^t\beta(\tau_0 + a)) = 1$ and let assumptions (C1)-(C4), (D1)-(D3) hold. Then the representation
\beq \label{betarep_bah}
\hat \beta(s) - \beta(s) = (\mu'(\beta(s)))^{-1} \Big(w_n(s) - \int_{[\tau_0,s)}\Big(\prodi_{(u,s]}\Big(I_{d} + M_v^tdH(v) \Big)\Big)^t M_u w_n(u)dH(u)\Big) + R_n(s)
\eeq
holds uniformly in $s \in [\tau_L,\tau_U]$ where $M_u = (\mu'(\beta(s)))^{-1}\tilde\mu'(\beta(u))$, $\prodi$ denotes the product-integral [see \cite{gilljoha1990}], and  for any $c_n \rightarrow \infty$ the remainder $R_n(s)$ satisfies
\[
\sup_{\tau \in [\tau_L,\tau_U]} \sqrt{n}\|R_n(\tau)\| = O_P(n^{1/2}b_n + n^{-\gamma/2} + \omega_{c_nn^{-1/2}}(\sqrt{n}\nu_n) + \omega_{c_nn^{-1/2}}(\sqrt{n}\tilde\nu_n))
\]
In particular, this implies
\bea
\sqrt n (\hat\beta(\cdot) - \beta(\cdot)) \Dkonv (\mu'(\beta(\cdot)))^{-1}\VVV_{\tau_0}(\cdot)
\eea
in the space $D([\tau_L,\tau_U])^{d}$ equipped with the supremum norm and ball sigma algebra [see \cite{pollard1984}]. Here $\VVV_{\tau_0}$ denotes centered Gaussian processes given by
\bean
\VVV_{\tau_0}(\tau) &=& \WWW(\tau) - \int_{[\tau_0,\tau)}\Big(\prodi_{(u,\tau]}\Big(I_{d} + M_v^tdH(v) \Big)\Big)^t M_u\WWW(u)dH(u).
\eean
\end{theo}

\medskip

The uniform Bahadur representation derived above has many potential applications. For example, it could be used to extend the L-statistic approach of \cite{kopor1987}, the rank tests of \cite{gujukopo1993}, or the confidence interval construction of \cite{zhoupor1996} to the setting of censored and/or dependent data. We conclude this section by discussing some interesting special cases and also possible extensions of the above result.

\medskip

\begin{rem}
{\rm
In the case of independent data, standard arguments from empirical process theory imply
\[
\omega_{c_nn^{-1/2}}(\sqrt{n}\nu_n) + \omega_{c_nn^{-1/2}}(\sqrt{n}\tilde\nu_n) = O_P(n^{-1/4}(c_n\log n)^{1/2}).
\]
Since $c_n$ can converge to infinity arbitrarily slowly, this shows that the remainder in (\ref{bahadur}) is of order $O_P(b_n+n^{-\gamma/2} + n^{-3/4}(\log n)^{1/2})$. In particular, for $\gamma \geq 1/2$ and $b_n = O(n^{-3/4})$ we obtain the same order as in the Bahadur representation of classical regression quantiles, see e.g. \cite{kopor1987}.
}
\end{rem}

\medskip

\begin{rem}
{\rm
If only conditions (D1) and (C1)-(C4) hold, the proofs of the result yield uniform consistency of the proposed quantile estimators. If the $o_p(1)$ in condition (D1) can be replaced by a rate $O_P(r_n)$ with $r_n$ tending to zero not faster then $n^{-1/2}$, it  is again possible to show that the censored regression quantiles converge uniformly with rate $O_P(b_n + r_n)$.
}
\end{rem}

\medskip

\begin{rem} {\rm
If there is no censoring we have $Y_i = X_i, \delta_i = 1, \ i=1,...,n$. In this case  $M_v= - I_d$ and thus for $0 < u \leq s < 1$
\[
\prodi_{(u,s]}\Big(I_{d} + (\M_v\tilde\M_v)^tdH(v) \Big) = \prodi_{(u,s]}(1-dH(v))I_d = \exp(H(s)-H(u))I_d = \frac{1-s}{1-u}I_d.
\]
In particular, in this case
\[
v_{\tau_0}(s) = w_n(s) - (1-s)\int_{[\tau_0,s)}\frac{w_n(u)}{(1-u)^2}du = w_n(\tau_0)\frac{1-s}{1-\tau_0} + \int_{[\tau_0,s)}\frac{1-s}{1-u}dw_n(u).
\]
After noting that $\ZZ^t\beta(\tau) = F_{Y|Z}^{-1}(\tau)$, and thus $I\{X_i \leq \ZZ^t\beta(u)\} = I\{F_{Y|Z}(X_i|\ZZ_i) \leq u\}$, straightforward but tedious calculations show that for $\delta_i \equiv 1$
\[
\int_{[0,s)}\frac{dw_n(u)}{1-u} = -\frac{1}{n}\sum_{i=1}^n \ZZ_i(I\{Y_i \leq \ZZ_i^t\beta(s)\} - s),
\]
which gives
\[
v_{\tau_0}(s) = -\frac{1}{n}\sum_{i=1}^n \ZZ_i(I\{Y_i \leq \ZZ_i^t\beta(s)\} - s).
\]
Thus the representation in (\ref{betarep}) corresponds to the Bahadur representation of regression quantiles in the completely uncensored case [see e.g. \cite{kopor1987}], and the proposed procedure is asymptotically equivalent to classical quantile regression.
}
\end{rem}

\section{Penalizing quantile processes} \label{genpen}

In this section we will discuss several aspects of penalization for quantile processes. For this purpose  we need   some additional notation and assumptions.
Let $\|\cdot\|$ denote the maximum norm in an Euclidean space. For a set $\mathcal{J} = \{j_1,...,j_J\} \subset \{1,...,d\}$ with $j_1<j_2<...<j_J$ define
\[
\beta^{(\mathcal{J})}=(\beta_j I \{j \in \mathcal{J} \})_{j=1,\dots,d}
\]
as the vector obtained from $\beta$, where components corresponding to indices $j \neq \mathcal{J}$ are set to zero. The vector $\bar \beta^{(\mathcal{J})}=(\beta_{j_1},...,\beta_{j_J})^t$ is defined as the vector of non-vanishing components of $\beta^{(\mathcal{J} )}$. Finally, introduce the matrix $\Pp_\mathcal{J}$ that corresponds to mapping coordinate $j_l$ to coordinate $l$ ($l=1,...,J$) and the remaining coordinates to $J+1,...,d$ (in increasing order). \\

Assume that the penalization in (\ref{intesteq.ava}) satisfies the following assumption (here $\Pp(A)$ denotes the power set of $A$)
\begin{enumerate}
\item[(P)] There exists a (set-valued) mapping $\chi: [\tau_L,\tau_U] \rightarrow \Pp(\{1,...,d\})$  such that
$\beta_k(\tau) = 0$ for all $k \in \chi(\tau)^C$, $\tau\in[\tau_L,\tau_U]$ and additionally
\bea \label{penas}
\sqrt{n} \Lambda_{n,0}  &:=& \sqrt{n}\inf_j \inf_{k \in \chi(\tau_j)^C} \frac{\lambda_n}{p_k(n,\tau_j)} \Pkonv \infty, \\
\Lambda_{n,1} &:=& \sup_j \sup_{k \in \chi(\tau_j)}\frac{\lambda_n}{p_k(n,\tau_j)} = o_P(1/\sqrt n). \label{penas1}
\eea
Moreover, there exist real numbers $\tau_L = \theta_1 < ... < \theta_K = \tau_U$ such that $\chi$ is constant on intervals of the form $[\theta_j,\theta_{j+1}), j=1,...,K-1$.
\end{enumerate}

A more detailed discussion of various penalizations satisfying condition (P) will be given in Sections \ref{detailsadapt} and \ref{detailsint}. In particular, in Section \ref{detailsadapt} we will provide conditions which guarantee that the adaptive lasso penalty   in (\ref{adappen}) fulfills (P)   and discuss what happens if those conditions fail. Alternative ways of choosing the penalty that do not suffer from the same problem and additionally allow to investigate the impact of covariates on multiple quantiles will be considered in Section \ref{detailsint}.\\
\\
For the results that follow, we need to strengthen assumption (D1) to
\begin{enumerate}
\item[(D1')] With the notation (\ref{nu}) we have
\[
\sup_{\bb \in \R^d} \|\nu_n(\bb)\| + \sup_{\bb \in \R^d} \|\tilde\nu_n(\bb)\| = O_P(n^{-1/2})
\]
\end{enumerate} 

Strengthening (D1) allows us to replace assumption (C4) by the weaker, and more realistic, version [note that for any $J \subset \{1,...,d\}$ we have $\lambda_{\min}(\E[(\bar\ZZ^{(J)})(\bar\ZZ^{(J)})^tf(\ZZ^t\bb|\ZZ)]) \geq \lambda_{\min}(\E[\ZZ\ZZ^tf(\ZZ^t\bb|\ZZ)])$ due to the special structure of the matrices].
\begin{enumerate}
\item[(C4')] We have for the map $\chi$ from condition (P)
\[
\inf_{\bb \in \B(\T,\eps)}\lambda_{\min}(\E[(\bar\ZZ^{(\chi(\tau))})(\bar\ZZ^{(\chi(\tau))})^tf(\ZZ^t\bb|\ZZ)]) =: \lambda_0 >0
\]
where $\lambda_{\min}(A)$ denotes the smallest eigenvalue of the matrix $A$.
\end{enumerate}

\begin{rem}{\rm
As discussed in Remark \ref{rem:vc}, the statement of (D1) can be viewed as a Glivenko-Cantelli type result for an empirical process indexed by a VC-subgraph class of functions. Similarly, (D1') follows if the same class of functions satisfies a Donsker type property. Results of this kind have for example been established for $\beta$-mixing data. More precisely, Corollary 2.1 in \cite{aryu1994} shows that (D1') holds as soon as the $\beta$-mixing coefficient $\beta_r$ satisfies $\beta_r = o(r^{-k})$ for some $k>1$.
}
\end{rem}

\begin{rem}{\rm
The results that follow continue to hold if we strengthen assumption (C4') to (C4) and replace (D1') by (D1). The details are omitted for the sake of brevity.}
\end{rem}
We now are ready to state our first  main result, which shows that under assumption (P) on the penalization, the estimate defined in (\ref{intesteq.ava}) enjoys the a kind of 'oracle' property in the sense of \cite{fanli2001}.
More precisely, with probability tending to one the coefficients outside the set $\chi(\tau)$ are set to zero uniformly in $\tau$
and the estimators of the remaining coefficients have the same asymptotic distribution as the estimators in the sub-model defined by $\chi(\tau)$.

\begin{theo} \label{satz1}
Assume that $\tau_0 = \tau_L >0$, that for some $a>0$ we have $P(C > \ZZ^t\beta(\tau_0 + a)) = 1$ and let assumptions (C1)-(C3),(C4'), (D1'), (D2)-(D3) and (P) hold. Then we have as $n \to \infty$
\bea \label{satz1a}
P(\sup_{\tau\in [\tau_L,\tau_U]} \sup_{k \in \chi(\tau)}|\hat\beta_k(\tau)| = 0) \rightarrow 1.
\eea
Moreover,
\beq \label{betarep}
\mu(\hat \beta(\tau)) - \mu(\beta(\tau)) = 
\M_{\tau,\chi} v_{\tau_L}(\tau) + o_P(1/\sqrt n)
\eeq
uniformly in  $\tau \in [\tau_L,\tau_U]$ where
\bean
v_{\tau}(s) &:=& w_n(s) - \int_{[\tau_0,s)}\Big(\prodi_{(u,s]}\Big(I_{d} + (\M_{v,\chi}\tilde\M_{v,\chi})^tdH(v) \Big)\Big)^t \tilde\M_{u,\chi} \M_{u,\chi} w_n(u)dH(u),
\eean
$\prodi$ denotes the product-integral [see \cite{gilljoha1990}], the matrices $\M_{\tau,\chi}, \tilde \M_{\tau,\chi}$ are defined by
\[
\M_{\tau,\chi} := \mu'(\beta(\tau))\Pp_{\chi(\tau)}^{-1}
\Big(
\begin{array}{cc}
M^{-1}_{\tau,\chi(\tau)} & 0\\
0 & 0
\end{array}
\Big)
\Pp_{\chi(\tau)}, \quad
\tilde\M_{\tau,\chi} := \tilde\mu'(\beta(\tau))\Pp_{\chi(\tau)}^{-1}
\Big(
\begin{array}{cc}
M^{-1}_{\tau,\chi(\tau)} & 0\\
0 & 0
\end{array}
\Big)
\Pp_{\chi(\tau)},
\]
and $M_{\tau,\chi(\tau)} := \E[(\bar\ZZ^{(\chi(\tau))})(\bar\ZZ^{(\chi(\tau))})^tf(\ZZ^t\beta(\tau)|\ZZ)]$. In particular, this implies
\bea
\sqrt n (\hat\beta(\cdot) - \beta(\cdot)) \Dkonv \Pp_{\chi(\cdot)}^{-1}
\Big(
\begin{array}{cc}
M^{-1}_{\cdot,\chi(\cdot)} & 0\\
0 & 0
\end{array}
\Big)
\Pp_{\chi(\cdot)}\M_{\cdot,\chi}\VVV_{\tau_0,\chi}(\cdot)
\label{satz1b}
\eea
in the space $(D[\tau_L,\tau_U])^{d}$ equipped with the supremum norm and ball sigma algebra [see \cite{pollard1984}]. Here $\VVV_{\tau_0}$ denotes a centered Gaussian process given by
\bean
\VVV_{\tau_0,\chi}(\tau) &=& \WWW(\tau) - \int_{[\tau_0,\tau)}\Big(\prodi_{(u,\tau]}\Big(I_{d} + (\M_{v,\chi}\tilde\M_{v,\chi})^tdH(v) \Big)\Big)^t \tilde\M_{u,\chi}\M_{u,\chi}\WWW(u)dH(u).
\eean
\end{theo}

\medskip

The asymptotic representation of the limiting process above is quite complicated. We now give a brief discussion of some special cases where it can be further simplified.

\medskip

\begin{rem} {\rm If there is no penalization, then $\chi(\tau) \equiv \{1,...,d\}$ and $\Pp_{\chi(\tau)}$ and $\M_{\tau,\chi}$ both are equal to the $d\times d$ identity matrix and $\tilde\M_{\tau,\chi} = \tilde\mu'(\beta(\tau))(\mu'(\beta(\tau)))^{-1}$. In this case, an analogue of Theorem \ref{satz1} is obtained from Theorem \ref{bahadur}, but without the rate on the remainder term.
If only the first $k<d$ components are important, i.e. if $\chi(\tau) \equiv \{1,...,k\}$ for $\tau \in [\tau_L,\tau_U]$, $\Pp_{\chi(\tau)}$ has a $k \times k$ identity matrix as the left upper block and the remaining entries are zero. The same holds for $\M_\tau$. Thus in this case the asymptotic distribution of the first $k$ components would be equal to the distribution in a smaller model where only those components are considered. This means that   the proposed procedure has a kind of 'oracle property'.
}
\end{rem}

\medskip

\begin{rem} {\rm
Under additional regularity assumptions, similar results can be derived for the version of the estimator starting with $\tau_0 = \tau_L = 0$ [see Remark \ref{rem:ph}]. The technical details are omitted for the sake of brevity.
}
\end{rem}

\subsection{Adaptive lasso penalization} \label{detailsadapt}

Recall the definition of the penalization in (\ref{adappen}) and assume that for some $\mathcal{J} \subset \{1,...,d\}$
\beq\label{noncross}
\inf_{k\in \mathcal{J}}\inf_{\tau \in [\tau_0,\tau_U]} |\beta_k(\tau)| > 0, \quad \sup_{k\in \mathcal{J}^C}\sup_{\tau \in [\tau_0,\tau_U]} |\beta_k(\tau)| = 0,
\eeq
then the following statement is correct.

\begin{cor} \label{cor1}
Assume that the conditions of Theorem \ref{satz1} are satisfied and that (\ref{noncross}) and
\beq \label{lamdban}
\sqrt n \lambda_n \rightarrow 0,~ n\lambda_n \rightarrow \infty
\eeq
hold. If the the preliminary estimator $\tilde\beta$ in (\ref{adappen}) is uniformly consistent with rate $O_P(1/\sqrt n)$ on the interval $[\tau_0,\tau_U]$, then the penalization (\ref{adappen}) satisfies (P) with $\chi(\tau) \equiv \mathcal{J}$. 
\end{cor}

The result shows that the adaptive lasso is $\sqrt{n}$ consistent under the assumption (\ref{noncross}).
It is of interest to investigate if a condition of this type is in fact necessary for the optimal rate of convergence.
The following result gives a partial answer to this question and  shows that the optimal rate cannot be achieved by the adaptive lasso defined in (\ref{adappen}) if some of the   coefficients
of the quantile regression change their sign or run into zero as $\tau$ varies.
More precisely we provide a lower bound on the uniform rate of convergence of the estimator
which turns out to be larger then $n^{-1/2}$ in quantile regions where coefficients are 'close' but not exactly equal to zero.
For a precise statement we define the sets [the dependence on $n$  is suppressed in the notation for the sake of brevity]
\bean
P_j &:=& \Big\{\tau\in [\tau_L,\tau_U]\Big|\frac{1}{n^{1/4}\kappa_n^{1/2}c_n} \leq |\beta_j(\tau)|\leq \frac{c_n}{\kappa_n} \Big\},\\
B_j &:=& \Big\{\tau\in [\tau_L,\tau_U]\Big| |\beta_j(\tau)| > \frac{c_n}{\kappa_n} \Big\},\\
S_{j} &:=& \Big\{\tau\in [\tau_L,\tau_U]\Big| \frac{1}{n^{1/4}\kappa_n^{1/2}c_n} > |\beta_j(\tau)| >0 \Big\},\\
V_j &:=& \Big\{\tau\in [\tau_L,\tau_U]\Big| \beta_j(\tau) = 0 \Big\}.
\eean

\begin{rem}{\rm
Basically, the sets defined above reflect the different kinds of asymptotic behavior of penalized estimators. The sets $B_j$ correspond to values of $\tau$ with $j$'th coefficients being 'large enough', such that they are not   affected by the penalization asymptotically. In contrast to that, coefficients $\beta_j(\tau)$ with $\tau \in S_j$ are 'too small' and will be set to zero with probability tending to one. In particular, this implies that the order of the largest elements in the set $\{|\beta_j(\tau)|: \tau \in S_j\}$ will give a lower bound for the uniform convergence rate of the penalized estimator. Finally, the set $P_j$ corresponds to 'intermediate' values that might be set to zero with positive probability.
}
\end{rem}

In order to state the next result, we need to make the following additional assumptions
\begin{enumerate}
\item[(C4*)] Define the map $\xi: [\tau_0,\tau_U] \rightarrow \Pp(\{1,...,d\})$ with $\xi(\tau) := \{j: |\beta_j(\tau)|\neq 0\}$.
Then
\[
\inf_{\bb \in \B(\T,\eps)}\lambda_{\min}(\E[(\bar\ZZ^{(\xi(\tau))})(\bar\ZZ^{(\xi(\tau))})^tf(\ZZ^t\bb|\ZZ)]) =: \lambda_0 >0
\]
where $\lambda_{\min}(A)$ denotes the smallest eigenvalue of the matrix $A$.
\item[(B1)] $\sqrt n \lambda_n = o(1), n\lambda_n \rightarrow \infty$, $\sqrt  n\kappa_n\lambda_n \rightarrow 1$
\item[(B2)] The set $P \cup S$ with $P := \cup_j P_j, S := \cup_j S_j$ is a finite union of intervals and its Lebesgue measure is bounded by $C_\gamma\Big(\frac{c_n}{\kappa_n} \Big)^{\gamma}$ for some positive constants $\gamma < \infty$ and a finite constant $C_\gamma$.
\item[(B3)] $c_n\rightarrow\infty$, $\lambda_n n^{3/4}\kappa_n^{1/2}c_n^{-1} \rightarrow \infty$, $n^{1/4}c_n^{\gamma+1}/\kappa_n^{\gamma+1/2} = o(1)$.
\item[(B4)] The preliminary estimator $\tilde \beta$ is uniformly consistent with rate $O_P(1/\sqrt n)$.
\end{enumerate}

\begin{rem}{\rm \label{rem:b1}
Assume that $\lambda_n \sim n^{-b}$ for some $b \in (1/2,1)$ and $c_n \sim \log(n)$ (it will later become apparent why choosing $c_n$ to converge to infinity slowly makes sense). Then $\kappa_n \sim n^{b-1/2}$, $\lambda_n n^{3/4}\kappa_n^{1/2} \sim n^{(1-b)/2}$ and $n^{1/4}/\kappa_n^{\gamma+1/2} \sim n^{(1+\gamma - b(\gamma+2))/2}$. Thus condition (B3) will hold as soon as $\frac{1}{2} \vee \frac{1+\gamma}{1+2\gamma} < b < 1$.
}
\end{rem}

\begin{rem}{\rm
Condition (B2) places a restriction on the behavior of the coefficients $\beta_j(\tau)$ in a neighborhood of $\{\tau|\beta_j(\tau)=0\}$. Essentially, it will hold if no coefficient approaches zero in a 'too smooth' way. If for example the function $\tau \mapsto \beta_j(\tau)$ is $k$ times continuously differentiable, (B2) will hold with $\gamma = 1/a$ where $a$ is the smallest number, such that the $a$'th derivative of $\beta_j(\tau)$  does not vanish at all points $\theta$ with $\beta_j(\theta) = 0$ for some $j$. In particular, in the case $\gamma = 1$ this property means that   $\beta(\tau)$ crosses zero with a positive slope. The results in Remark \ref{rem:b1} show that $\lambda_n \sim n^{-b}$ for any $b \in (1/2,1)$ is allowed when $c_n = \log n$. If $\beta(\tau)$ runs into zero more smoothly, which corresponds to $\gamma<1$, the conditions on the regularizing parameter $\lambda_n$ become stricter since now only $\frac{1}{2} \vee \frac{1+\gamma}{1+2\gamma} < b < 1$ is allowed.
}
\end{rem}

\begin{theo}  \label{lem:unifrates}
Assume that conditions (C1)-(C3), (C4*), (D1'), (D2)-(D3), (B1)-(B4) hold. Then adaptive lasso estimator obtained form the penalization (\ref{adappen}) satisfies
\beq \label{alarate}
\sup_{\tau \in [\tau_L,\tau_U]} \|\hat \beta(u) - \beta(u)\| = O_P\Big(\frac{c_n}{\kappa_n^{1/2}n^{1/4}}\Big).
\eeq
Moreover, for any fixed $I \subset [\tau_L,\tau_U] \backslash (S \cup P)$
\beq \label{alawk}
\sqrt{n}(\hat \beta(\cdot) - \beta(\cdot)) \rightarrow \Pp_{\xi(\tau)}^{-1}
\Big(
\begin{array}{cc}
M^{-1}_{\tau,\xi(\tau)} & 0\\
0 & 0
\end{array}
\Big)
\Pp_{\xi(\tau)}\VVV_{\tau_0,\xi}(\cdot)
\eeq
in the space $D(I)^d$ where the process $\VVV_{\tau_0,\xi}$ is defined   in Theorem \ref{satz1} and
\beq \label{alazero}
\Pe(\sup_{j=1,...,d}\ \sup_{\tau\in S_{j}\cup V_j\cap[\tau_L,\tau_U]}|\hat\beta_j(\tau)|=0) \rightarrow 1.
\eeq
\end{theo}

Note that the assertion (\ref{alazero}) implies that the uniform rate of $\hat \beta$ is bounded from below by $n^{-1/4}\kappa_n^{-1/2}c_n^{-1}$ as soon as the set $S\cup P$ is not empty. Since $c_n$ is allowed to converge to infinity arbitrarily slow, we obtain the lower bound $O(n^{-1/4}\kappa_n^{-1/2}) = O(\lambda_n^{1/2})$, which depends on $\lambda_n$ and is always slower then $1/\sqrt n$. We will demonstrate in Section \ref{simul} by means of a simulation study that this inferior property of the adaptive lasso can also be observed for realistic sample sizes.

\begin{rem}{\rm
Theorem \ref{lem:unifrates} also contains a positive, and at the first glance probably surprising, result. Since the procedure used to compute the estimators is iterative, one might expect that a non-optimal convergence rate of the estimator at one value of $\tau$ should yield the same lower bound for all subsequent quantile estimators. However, the above results imply that this is not always the case. The intuitive reason for this phenomenon is the following: the estimators $\hat\beta(\tau)$ only enter the subsequent estimating equation inside an integral, see equation (\ref{intesteq.ava}). Thus, when the rate is not optimal on a sufficiently small set of values $\tau$, the overall impact of a non-optimal rate might still be small. In particular, this is the case under conditions (B2)-(B4).
}
\end{rem}

\begin{rem}{\rm
The results in the above Theorem are related to the findings of \cite{potle2009} which demonstrate that penalized estimators do not have optimal convergence rates uniformly over the parameter space. This also suggests that using other point-wise penalties such as for example SCAD will not solve the problems encountered by the adaptive lasso. Instead, using information from other quantiles is necessary.  
}
\end{rem}

\subsection{Average penalization} \label{detailsint}

\noindent
As we have seen in the last section, the traditional  way of implementing the adaptive lasso will yield sub-optimal rates of convergence if some coefficients cross zero. Moreover, this method will perform a 'point-wise' model selection with respect to quantiles- a property, which might not always be desirable. Rather, keeping the same model for certain ranges of quantiles such as for example $\tau \in [.4,.6]$, or even for the whole range, might often be preferable. In order to implement such an approach, and to obtain a quantile process which converges at the optimal rate, we introduce a new kind of adaptive penalization which has - to the best of our knowledge - not been considered in the literature so far . More precisely, denote by $\T_1,...,\T_K$ a fixed, disjoint partition of $[\tau_0,\tau_U]$ and define
\bea\label{pen.ava}
p_k^{int}(n,\tau) &:=& \sum_{j=1}^K I\{\tau \in \T_j\} \int_{\T_j} | \tilde \beta_k(t)| h(t) dt, \quad k=1,...,d
\\
p_k^{max}(n,\tau) &:=& \sum_{j=1}^K I\{\tau \in \T_j\} \sup_{t\in\T_j} |\tilde \beta_k(t)| , \quad k=1,...,d. \label{pen.ava2}
\eea
Here, $\tilde \beta$ is a preliminary estimator which converges uniformly with rate $O_P(1/\sqrt n)$ on the interval $[\tau_0,\tau_U]$, and $h$ is a strictly positive, uniformly bounded weight function integrating to one. In the following discussion we call this method \textit{average adaptive lasso.}

\begin{rem}{ \rm
The above idea can be generalized to the setting where the researcher wants to include a whole set of predictors, say $(Z_k)_{k \in S}$, in the analysis if at least one of those predictors is important. This can be done by setting
\[
p_k(n,\tau_j) := \max_{m \in S} \sum_{j=1}^K I\{\tau \in \T_j\} \sup_{t\in\T_j} |\tilde \beta_m(t)| , \quad k \in S.
\]
}
\end{rem}

\begin{rem}{\rm
In the context of uncensored quantile regression, \cite{zouyuan2008} recently proposed to simultaneously penalize a collection of estimators for different quantiles in order to select the same group of predictors for different values of the quantile. While such an approach is extremely interesting, it seems hard to implement in the present situation. The reason is that the minimization problem (\ref{intesteq.ava}) is solved in an iterative fashion, and dealing with a penalty that affects all quantiles at the same time thus is problematic.
}
\end{rem}

The following result follows

\begin{lemma}
\label{cor2}
Assume that there exist sets $\mathcal{J}_1,...,\mathcal{J}_K \subset \{1,...,d\}$ such that
\beq \label{intcond}
\inf_{j=1,...,K}\inf_{k\in \mathcal{J}_j}\sup_{\tau \in \T_j} |\beta_k(\tau)| > 0, \quad \sup_{j=1,...,K}\sup_{k\in \mathcal{J}_j^C}\sup_{\tau \in \T_j} |\beta_k(\tau)| = 0,
\eeq
and (\ref{lamdban}) hold. If the the preliminary estimator $\tilde\beta$ is uniformly consistent with rate $O_P(1/\sqrt n)$ on the interval $[\tau_{L},\tau_{U}]$ then the average penalties defined in (\ref{pen.ava}) and (\ref{pen.ava2}) satisfy (P) with $\chi(\tau) = \mathcal{J}_j$ for $\tau \in \T_j$.
\end{lemma}

The above results imply that the problems encountered by the traditional  application of adaptive lasso when coefficients cross zero can be avoided if average penalization is used. Another consequence of such an approach is that predictors which are important for some quantile $\tau \in \T_k$ will be included in the analysis for all quantiles in $\T_k$. At the same time, covariates that have no impact for any $\tau \in \T_k$ can still be excluded from the analysis. Finally, by taking $\T_1 = [\tau_0,\tau_U]$ it is possible to achieve that all covariates that are important at some quantile in the range of interest will be used for all $\tau \in [\tau_0,\tau_U]$. As a consequence, average penalization is a highly flexible method that can easily be adapted to the situation at hand.

\section{Simulation study} \label{simul}

In order to study the finite-sample properties of the proposed procedures  we conducted a small simulation study. An important practical question is the selection of the regularizing parameter $\lambda_n$. In our simulations, we used an adapted version of $K$-fold cross validation  which accounts for the presence of censoring by using a weighted objective function. More precisely we proceeded is two steps. In the first step, weights were estimated as follows
\begin{enumerate}
\item Compute an unpenalized estimator based on all data, denote this estimator by $\hat\bb$.
\item For each grid point $\tau$, following \cite{portnoy2003} define weights $\hat w_j(\tau)$ through
\[
\hat w_j(\tau) := \delta_j + (1-\delta_j)\Big(I\{X_j > \ZZ_i^t \hat\bb(\tau)\} + I\{X_j \leq \ZZ_i^t \hat\bb(\tau)\}\frac{\tau-r_j}{1-r_j} \Big)
\]
Here, $r_j$ denotes the value of $\tau$ at that the observation $X_j$ is 'crossed', that is
\[
r_j :=
\left\{
\begin{array}{ccr}
1 &\mbox{if}& X_j > \ZZ_i^t \hat\bb(\tau_U)
\\
\inf\{\tau_k| \ZZ_i^t \hat\bb(\tau_{k-1}) < X_j \leq \ZZ_i^t \hat\bb(\tau_k)\} &\mbox{if}& \delta_j = 0, X_j \leq \ZZ_i^t \hat\bb(\tau_U)
\\
0 &\mbox{if}& \delta_j = 1, X_j \leq \ZZ_i^t \hat\bb(\tau_U)
\end{array}
\right.
\]
\end{enumerate}
Note that \cite{portnoy2003} used the weights $\hat w_j(\tau)$ to define a weighted minimization problem to account for censoring. The basic idea corresponds to the well-known interpretation of the classical Kaplan-Meier estimator as an iterative redistribution of mass corresponding to censored observations to the right. After obtaining preliminary estimators of the weights, the second step was to select $\lambda$ as the minimizer of the function $CV(\lambda)$ which was computed as follows.

\begin{enumerate}
\item Randomly divide the data into $K$ blocks of equal size. Denote the corresponding sets of indexes by $J_1,...,J_K$.
\item For $k=1,...,K$, compute estimators $\hat\bb^{(J_k,\lambda)}$ based on the data $(\ZZ_i,X_i,\delta_i)_{i \in \{1,...,n\}\backslash J_k}$ and penalization level $\lambda$.
\item Compute
\[
CV(\lambda) := \sum_{k=1}^K \sum_{j \in  J_k}\sum_{i =1}^{N_\tau(n)} \Big(\hat w_j(\tau_i)\rho_{\tau_i}(X_j - \ZZ_j^t\hat\bb^{(J_k,\lambda)})
+ (1-\hat w_j(\tau_i))\rho_{\tau_i}(X^\infty - \ZZ_j^t\hat\bb^{(J_k,\lambda)})\Big)
\]
where $X^\infty$ denotes some sufficiently large number (we chose $10^3$ in the simulations). Select the penalty parameter $\lambda$ as the minimizer of $CV(\lambda)$ among a set of candidate parameters.
\end{enumerate}

The basic idea behind the above procedure is that the weights $\hat w_i$ are consistent 'estimators' of the random quantities
\[
w_i(\tau_j) = \delta_i + (1-\delta_i)\Big(I\{X_i > F_T^{-1}(\tau_j|\ZZ_i)\} + I\{X_i \leq F_T^{-1}(\tau_j|\ZZ_i)\}\frac{\tau-F_T(X_i|\ZZ_i)}{1-F_T(X_i|\ZZ_i)} \Big),
\]
and that the minimizer of the weighted sum
\[
\sum_{j =1}^{n} \Big(w_j(\tau)\rho_{\tau_i}(X_j - \ZZ_j^t\bb)
+ (1-w_j(\tau))\rho_{\tau_i}(X^\infty - \ZZ_j^t\bb)\Big)
\]
is a consistent estimator of $\beta(\tau)$. See \cite{portnoy2003} for a more detailed discussion.

\begin{rem} {\rm
At the first glance, it might seem that by redistributing mass to $X^\infty$ we would give higher quantiles more importance since the corresponding quantile curves have crossed more observations. However, while it is true that the total value of the sum
\[
\sum_{j \in  J_k}\sum_{i =1}^{N_\tau} \Big(\hat w_j(\tau_i)\rho_{\tau_i}(X_j - \ZZ_j^t\hat\bb^{(J_k,\lambda)})
+ (1-\hat w_j(\tau_i))\rho_{\tau_i}(X^\infty - \ZZ_j^t\hat\bb^{(J_k,\lambda)})\Big)
\]
will be larger for higher quantiles, the magnitude of changes induced by perturbations of $\lambda$ will in fact be of the same order across quantiles. In a certain sense, this corresponds to the invariance of regression quantiles to moving around extreme observations.}
\end{rem}

We considered two models. In the first model,     we generated data from
\[
\mbox{(model 1)} \quad
\left\{
\begin{array}{ccc}
T_i &=& (Z_{i,2},...,Z_{i,10})\tilde\bb + .75U_i
\\
C_i &=& (Z_{i,2},...,Z_{i,10})\tilde\bb + .75V_i
\end{array}
\right.
\]
where $\tilde\bb = (.5,1,1.5,2,0,0,0,0,0)^t$, $Z_{i,2},...,Z_{i,10}$ are independent $\mathcal{U}[0,1]$ distributed random variables and $U_i,V_i$ are independent $\mathcal{N}(0,1)$. The amount of censoring is roughly $25\%$. In this model, all coefficients are bounded away from zero and so the local adaptive lasso as well as the average  penalization methods share the same $n^{-1/2}$ convergence rates. We estimated the quantile process based on the grid $\tau_L = .15, \tau_U = .7$ with steps of size $.01$. Our findings are summarized in Table \ref{mod1}, which shows the integrated [over the quantile grid] mean squared error (IMSE) and the probabilities of setting coefficients to $0$ for the two estimates obtained by the different penalization techniques. All reported results are based on 500 simulation runs and $K=5$ in the cross validation. Overall, both estimators behave reasonably well. The average penalization method is always at least as good as the local penalization method. It has a systematically higher probability of setting zero components to zero and a systematically lower IMSE for estimating the intercept and the coefficient $\beta_2$.

\begin{table}[ht]
\center{
\begin{tabular}[c]{|c|c|cccccccccc|}
\hline
$n$& method & $\beta_{1}$ & $\beta_{2}$ & $\beta_{3}$ & $\beta_{4}$ & $\beta_{5}$ & $p_2$ & $p_3$ & $p_4$ & $p_5$ & $p_0$\\
\hline
100 & local &33.0& 16.8& 19.6& 17.3& 15.8&  40.6&  5.3&  0.1&  0.0& 71.7\\
&average &31.1& 16.2& 18.4& 16.6& 15.4&  40.9&  3.4&  0.0&  0.0& 75.8\\
\hline
250 &local&29.1& 20.8& 17.7& 14.9& 13.2&  15.2&  0.1&  0.0&  0.0& 72.2\\
&average&27.3& 19.9& 17.1& 14.8& 13.2&  13.4&  0.0&  0.0&  0.0& 77.1\\
\hline
500 &local&21.8& 19.9& 13.2& 13.3& 13.5&   3.6&  0.0&  0.0&  0.0& 76.7\\
&average&20.0& 18.4& 13.1& 13.2& 13.4&   2.2&  0.0&  0.0&  0.0& 80.8\\
\hline
1000 &local&20.8& 17.2& 13.3& 12.5& 12.3&   0.1&  0.0&  0.0&  0.0& 80.6\\
&average&19.3& 16.0& 13.0& 12.5& 12.3&   0.0&  0.0&  0.0&  0.0& 84.7\\
\hline
\end{tabular}}
\caption{\label{mod1} \emph{Results for model 1. Columns 1-5 show $n*IMSE(\beta_j), j=1,\dots,5$, where $\beta_1$ corresponds to the intercept. Columns 6-9 show the probabilities $p_j$ of  setting the coefficient $\beta_j$ to zero ($j=2,\dots,5$) averaged over all quantiles on the grid. Column 10 shows the average probability $p_0$ of setting coefficients $\beta_6-\beta_{10}$ to zero. Rows with label 'local' correspond to (local) adaptive lasso, rows with label 'average' correspond to average adaptive lasso.}}
\end{table}

\begin{center}
\begin{figure}
 \includegraphics[width=0.45\textwidth]{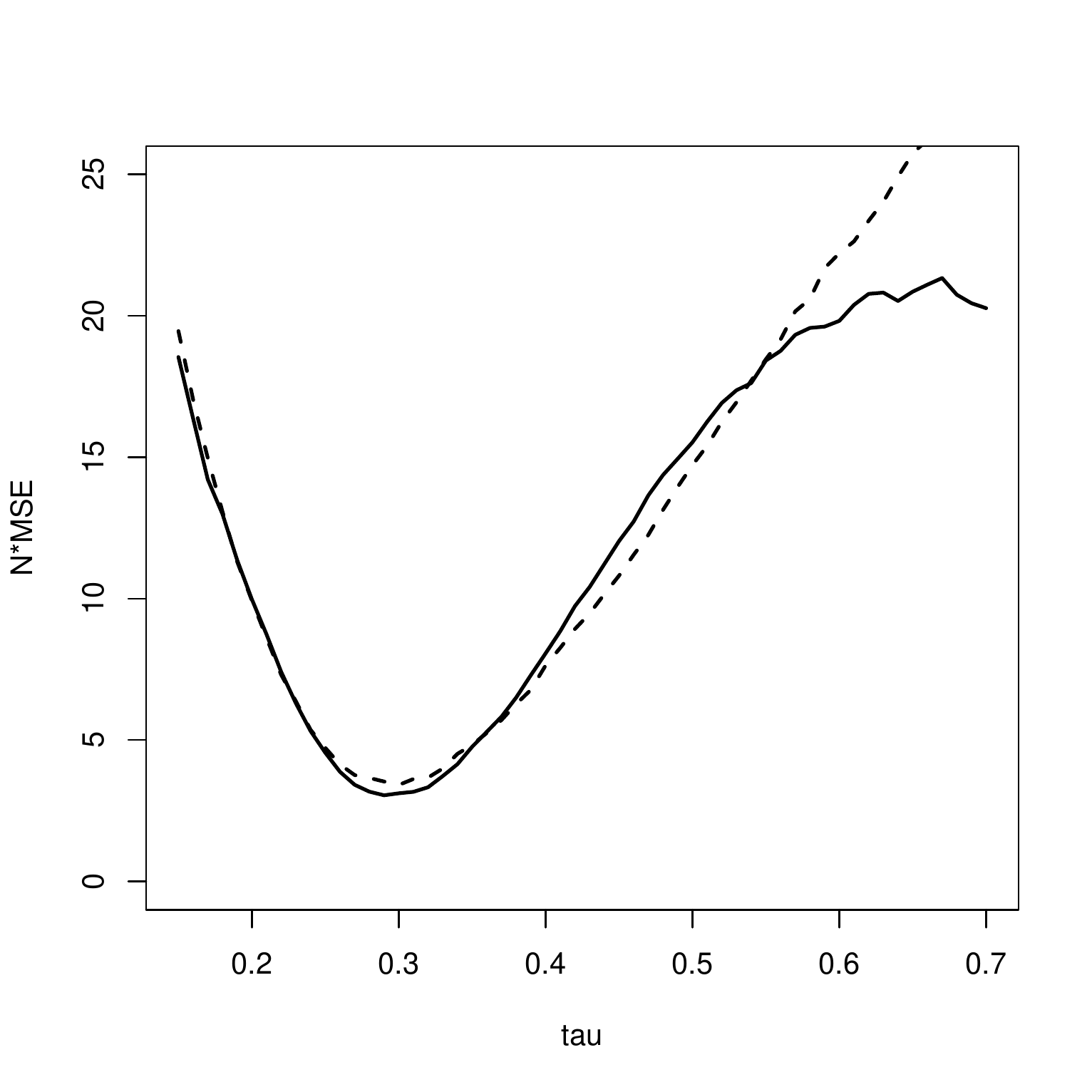}
 \includegraphics[width=0.45\textwidth]{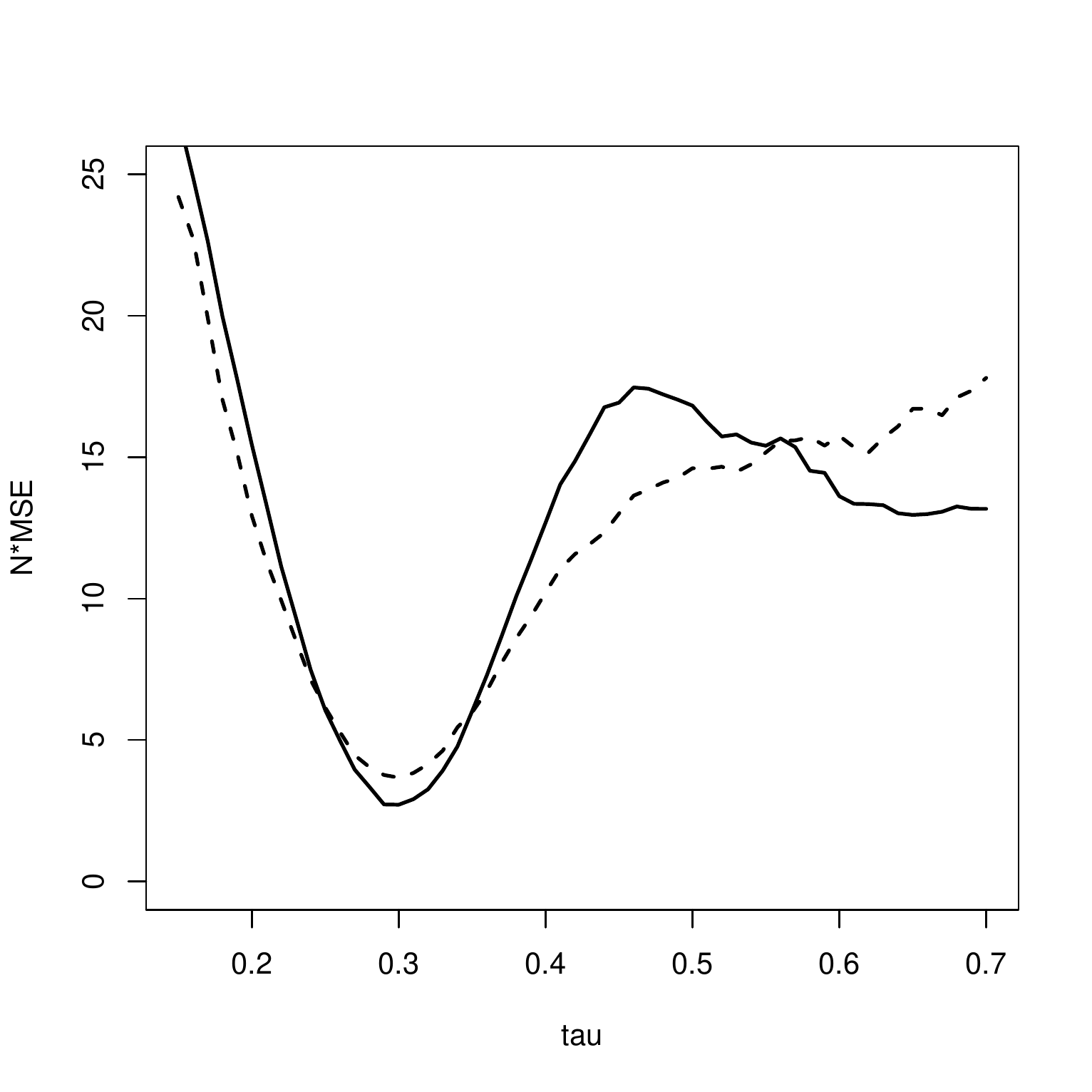}
 \includegraphics[width=0.45\textwidth]{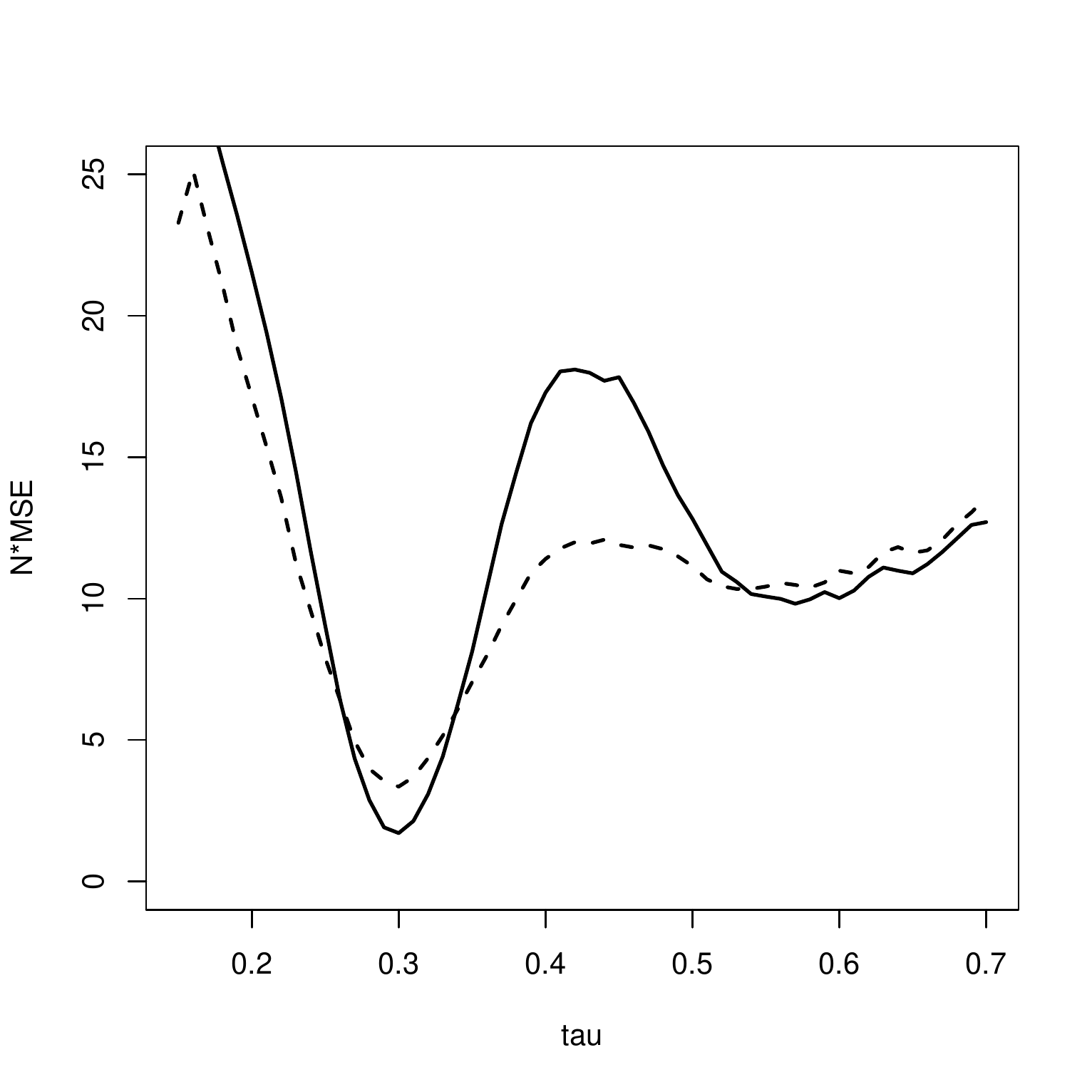}
 \includegraphics[width=0.45\textwidth]{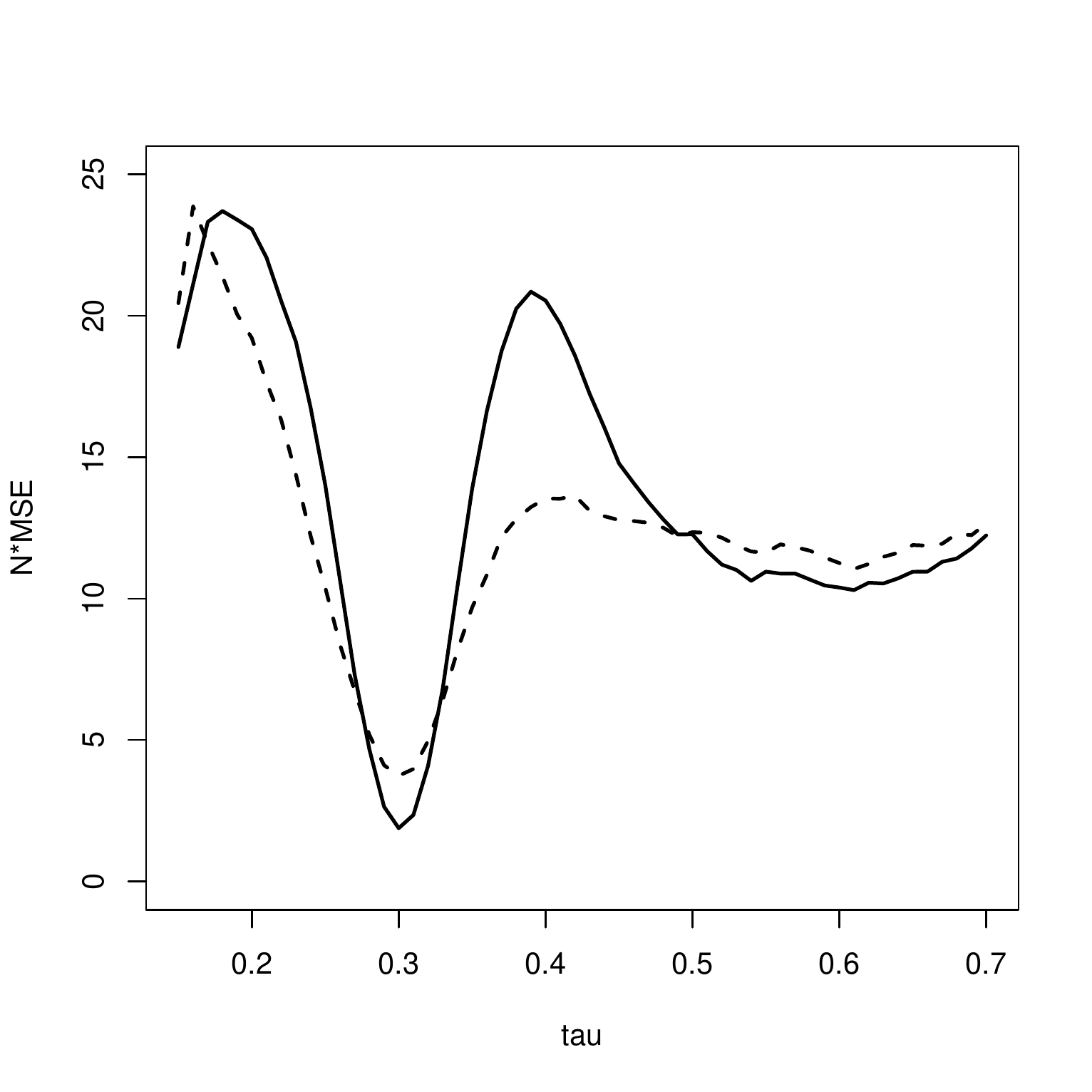}
\caption{\label{fig1} \emph{n*MSE of the estimate for the coefficient $\beta_7$ as a function of the quantile for sample sizes $n=50$ (upper left), $n=100$ (upper right), $n=250$ (lower left) and $n=1000$ (lower right). Solid line: local penalization. Dashed line: average penalization.}}
\end{figure}
\end{center}

The second model was of the form
\[
\mbox{(model 2)} \quad
\left\{
\begin{array}{ccc}
T_i &=& (Z_{i,2},...,Z_{i,6})\tilde\bb + Z_{i,7}(U_i - q)
\\
C_i &=& (Z_{i,2},...,Z_{i,6})\tilde\bb + 1.5 + V_i
\end{array}
\right.
\]
where $q$ denotes the $30\%$-quantile of a standard normal random variable, $Z_{i,2},...,Z_{i,7}$ are independent, $.2+\mathcal{U}[0,1]$-distributed random variables, $U_i$, $V_i$ are independent $\mathcal{N}(0,1)$ distributed, and $\tilde\bb = (2,2,0,0,0)$. The amount of censoring is roughly $20\%$.  We have calculated the quantile regression estimate for the model 
$$Q_\tau(T_i | \ZZ_i)= \beta_1(\tau) + \sum^7_{j=2}\beta_j(\tau) Z_{i,j}.$$
In this model, the coefficient corresponding to $Z_{i,7}$ crosses zero for $\tau = 0.3$. From an asymptotic point of view the estimator based on point-wise penalization should thus have a slower rate of convergence in a neighborhood of $\tau = 0.3$. First, consider the results in Table \ref{mod2} for the IMSE and the probabilities of setting coefficients to $0$. We observe the same slight but systematic advantages for the average penalization method with respect to model selection properties and integrated MSE. Note that this is consistent with the theory since the range of quantiles where the local penalization has a slower rate of convergence is shrinking with $n$. Plotting the MSE of the estimator $\hat\beta_7$ as a function of $\tau$ reveals a rather different picture [see Figure \ref{fig1}]. Here, the suboptimal rate of convergence of the local
 penalization and the clear asymptotic superiority of the average penalization becomes apparent.
\begin{table}[ht]
\center{
\begin{tabular}[c]{|c|c|cccccccc|}
\hline
$n$& method & $\beta_{1}$ & $\beta_{2}$ & $\beta_{3}$ & $\beta_{7}$ &  $p_2$ & $p_3$ & $p_7$ & $p_0$\\
\hline
100 & local &17.4&  9.3&  9.5& 12.5&  0.0&  0.0& 44.8& 76.3\\
 & average &17.0&  9.1&  9.6& 13.4&  0.0&  0.0& 39.4& 79.7\\
\hline
250 & local &14.6&  7.7&  8.1& 12.7&  0.0&  0.0& 31.8& 79.4\\
 & average &14.0&  7.7&  8.1& 12.3&  0.0&  0.0& 19.8& 82.1\\
\hline
500 & local &14.0  &7.9  &7.9 &12.8  &0.0  &0.0 &23.9& 82.3\\
 & average &12.8  &7.9  &7.9 &11.2  &0.0  &0.0 &11.9 &87.0\\
\hline
1000 & local &13.8&  7.4&  7.2& 13.5&  0.0&  0.0& 17.4& 83.5\\
 & average &12.6&  7.4&  7.2& 12.3&  0.0&  0.0&  8.1& 90.7 \\
\hline
\end{tabular}}
\caption{\label{mod2} \emph{Results for model 2. Columns 1-4 show  $n*IMSE(\beta_j), j=1,2,3,7$, where $\beta_1$ corresponds to the intercept. Columns 5-7 show the probabilities $p_j$ of  setting the coefficient $\beta_j$ to zero ($j=2,3,7$) averaged over all quantiles on the grid. Column 8 shows the average probability $p_0$ of setting coefficients $\beta_4-\beta_6$ to zero. Rows with label 'local' correspond to (local) adaptive lasso, rows with label 'average' correspond to average adaptive lasso.}}
\end{table}

\newpage

\section{Appendix: proofs} \label{proofs}

At the beginning of the proofs, we give a brief overview of the main results. Several auxiliary results are proved in Section
\ref{aux}. A first key result  here is Lemma \ref{lem:bednull.gen} which provides some general bounds for $\mu_k(\hat\beta(\tau_j)) - \mu_k(\beta(\tau_j))$. Moreover, conditions that describe when coefficients $\hat\beta_k$ are set to zero are derived. Lemma \ref{lem:bednull.gen} will play a major role in the proof of the subsequent results. Lemma \ref{lem:dmuphi} shows that $\sqrt n(\mu(\hat\beta(\cdot)) -  \mu(\beta(\cdot)))$ is uniformly close to $\sqrt n(\mu(\hat\beta(\cdot)) - \phi_n(\tau_j))$, which in turn is obtained as the solution of an iterative equation. Thus the asymptotic distribution of the two aforementioned quantities coincide. We will then proceed in Lemma \ref{lem:phi} to derive an explicit, i.e. non-iterative, representation for the quantity $\sqrt n(\mu(\hat\beta(\cdot)) - \phi_n(\tau_j))$.  This will yield a Bahadur representation of the process $\sqrt n(\mu(\hat\beta(\cdot)) - \mu(\beta(\cdot))$, which in turn is the main ingredient for establishing the representation for $\sqrt n (\hat\beta(\cdot) - \beta(\cdot))$. Since the proofs of the results in Sections \ref{dependent} and \ref{genpen} are similar, we only give detailed arguments for the results in Section \ref{genpen} [which are more complicated] and briefly mention the differences where necessary.

\subsection{Preliminaries} \label{aux}
We begin by stating some useful technical facts and introducing some notation that will be used throughout the following proofs.
\begin{rem} {\rm \label{rem.mu} ~~~\\
(1) Under condition (C3) it follows that, for any $\bb_1,\bb_2 \in \B(\T,\eps)$, $ \|\mu(\bb_1)-\mu(\bb_2) \| \leq C_2\|\bb_1-\bb_2\| $ with $C_2 := dC_Z^2K_f$ and
$ \|\tilde\mu(\bb_1)-\tilde\mu(\bb_2)\| \leq C_3\|\bb_1-\bb_2\| $ with $C_3 := dC_Z^2 \tilde K_f$.\\
\\
(2) Condition (C4') implies the inequality
\[
\|\mu^{(\chi(\tau))}(\bb_1^{(\chi(\tau))}) - \mu^{(\chi(\tau))}(\bb_2^{(\chi(\tau))})\| \geq \frac{\lambda_0}{|\chi(\tau)|}\ell_\tau(\bb_1^{(\chi(\tau))},\bb_2^{(\chi(\tau))}) \|\bb_1^{(\chi(\tau))}-\bb_2^{(\chi(\tau))} \|.
\]
where $\ell_\tau(\bb_1^{(\chi(\tau))},\bb_2^{(\chi(\tau))}) := \lambda(\{\gamma \in [0,1]: \|\gamma\bb_1^{(\chi(\tau))}+(1-\gamma)\bb_2^{(\chi(\tau))} - \beta(\tau) \|\leq\eps \})$ and $\lambda$ denotes the Lebesgue measure (a sketch of the proof is below). In particular, the above equation implies that for all $\bb$ with $\|\bb^{(\chi(\tau))} - \beta(\tau)\| \leq \frac{\eps}{C_1 \vee 1}$ with $C_1 := d/\lambda_0$ it holds that
\[
\|\bb^{(\chi(\tau))} - \beta(\tau)\| \leq \frac{1}{C_1}\|\mu^{(\chi(\tau))}(\bb_1^{(\chi(\tau))}) - \mu^{(\chi(\tau))}(\beta(\tau))\|.
\]
For a proof of the inequality above, note that
\bean
&&|J| \|\bar \bb_1^{(J)} - \bar \bb_2^{(J)}\|\|\bar \mu^{(J)}(\bb_1^{(J)}) - \bar \mu^{(J)}(\bb_2^{(J)})\| \geq (\bar \bb_1^{(J)} - \bar \bb_2^{(J)})^t (\bar\mu^{(J)}(\bb_1^{(J)}) - \bar\mu^{(J)}(\bb_2^{(J)}))
\\
&& ~~~~  = \E[(\bar \ZZ^{(J)})^t(\bar \bb_1^{(J)} - \bar \bb_2^{(J)})(F(\ZZ^t \bb_1^{(J)}|\ZZ) - F(\ZZ^t \bb_2^{(J)}|\ZZ))]
\\
&& ~~~~  = \E\Big[(\bar \ZZ^{(J)})^t(\bar\bb_1^{(J)} - \bar \bb_2^{(J)})(\bar \ZZ^{(J)})^t(\bar \bb_1^{(J)} - \bar \bb_2^{(J)}) \int_0^1 f(\ZZ^t(\gamma\bb_1^{(J)}+(1-\gamma)\bb_2^{(J)})|\ZZ) d\gamma\Big]
\\
&& ~~~~   = \int_0^1 (\bar \bb_1^{(J)} - \bar \bb_2^{(J)})^t \E\Big[ (\bar \ZZ^{(J)}) (\bar \ZZ^{(J)})^t f(\ZZ^t(\gamma\bb_1^{(J)}+(1-\gamma)\bb_2^{(J)})|\ZZ) \Big] (\bar \bb_1^{(J)} - \bar \bb_2^{(J)}) d\gamma
\eean
(3) For $\|\bb^{(\chi(\tau))} - \beta(\tau)\|\leq \eps$ we have
\begin{eqnarray*}
\mu(\bb^{(\chi(\tau))}) - \mu(\beta(\tau)) &=& \M_\tau\Big(\bar\mu^{(\chi(\tau))}(\bb^{(\chi(\tau))}) - \bar\mu^{(\chi(\tau))}(\beta(\tau))\Big) + D_\tau(\bb),
\\
\tilde\mu(\bb^{(\chi(\tau))}) - \tilde\mu(\beta(\tau)) &=& \tilde\M_\tau\Big(\bar\mu^{(\chi(\tau))}(\bb^{(\chi(\tau))}) - \bar\mu^{(\chi(\tau))}(\beta(\tau))\Big) + \tilde D_\tau(\bb)
\end{eqnarray*}
where $\sup_\tau \|D_\tau(\bb)\| = O(\|\bb - \beta(\tau)\|^\gamma), \sup_\tau \|\tilde D_\tau(\bb)\|=O(\|\bb - \beta(\tau)\|^\gamma)$. Introduce the notation 
\beq\label{def.d}
\D(a) := \sup_{\tau \in [\tau_L,\tau_U]}\sup_{\|\bb - \beta(\tau)\|\leq a}\|D_\tau(\bb)\|,\quad \tilde \D(a) := \sup_{\tau \in [\tau_L,\tau_U]}\sup_{\|\bb - \beta(\tau)\|\leq a}\|\tilde D_\tau(\bb)\|.
\eeq
(4) Assumptions (C2)-(C4) imply the existence of finite constants $C_5,\tilde C_5$ such that for any $\|\bb^{(\chi(\tau))} - \beta(\tau)\| \leq \eps$ we have
\begin{eqnarray} \label{muineq1}
\|\mu(\bb^{\chi(\tau)}) - \mu(\beta(\tau))\| &\leq& C_5\|\bar\mu^{(\chi(\tau))}(\bb^{(\chi(\tau))}) - \bar\mu^{(\chi(\tau))}(\beta(\tau))\|,
\\ \label{mutineq1}
\|\tilde\mu(\bb^{\chi(\tau)}) - \tilde\mu(\beta(\tau))\| &\leq& \tilde C_5 \|\bar\mu^{(\chi(\tau))}(\bb^{(\chi(\tau))}) - \bar\mu^{(\chi(\tau))}(\beta(\tau))\|.
\end{eqnarray}
}
\end{rem}

\newpage

\begin{lemma}\label{lem:bednull.gen}
Let $J \subset \{1,...,d\}$, $\vec 0 \in \R^{d-|J|}$ and consider the problem of minimizing $H_j(\Pp_{J}^{-1}(\hh^t,\vec 0^t)^t)$ with respect to $\hh\in\R^{|J|}$. Denote the generalized solution of this minimization problem by $\hat \hh(\tau_j)$ and set $\check \beta=\Pp_{J}^{-1}(\hat \hh(\tau_j)^t,\vec{0}^t)^t$. Then 
\bean
&&\Big|\mu_k(\check\beta) - \mu_k(\beta^{(J)}(\tau_j)) + \nu_{n,k}(\check\beta) - \int_{[\tau_0,\tau_j)} \tilde\nu_{n,k}(\hat\beta(u))dH(u) - \int_{[\tau_0,\tau_j)} \tilde \mu_k(\hat\beta(u)) - \tilde \mu_k(\beta(u)) dH(u)
\\
&&\quad - \frac{\tau_0}{n}\sum_{i=1}^n (\ZZ_{i,k} - \E\ZZ_{i,k})\Big|
\\
&\leq& \frac{\lambda_n}{p_k(n,\tau_j)} + \frac{C_Z}{n} + \|\mu^{(J)}(\beta(\tau_j)) - \mu^{(J)}(\beta^{(J)}(\tau_j))\|
\eean
for $k\in J$. \\
\\
Now, let conditions (P), (C1)-(C3) and (C4') hold and additionally assume that for some $J\subset \chi(\tau_j)$
\bea
\label{bednull.ava}
&& \quad \alpha_1 + \alpha_2 + \sup_{k \in J}\frac{\lambda_n}{p_k(n,\tau_j)} + \frac{C_Z}{n}
+ C_2 \sup_{k \in J^C} |\beta_k(\tau_j))| \leq \frac{\eps - \sup_{k \in J^C} |\beta_k(\tau_j))|}{C_1\vee 1}
\\
\label{bednull2.ava}
&& \sup_{k \in J}\frac{C_5\lambda_n}{p_k(n,\tau_j)}+ (C_5+1)\Big(\frac{2C_Z}{n} + \alpha_1 + \alpha_2 + C_2 \sup_{k \in J^C} |\beta_k(\tau_j))| \Big)\leq \inf_{k\in J^C} \frac{\lambda_n}{p_{k}(n,\tau_j)} \quad\quad\quad\quad
\eea
where
\bean
\sup_{\bb\in\R^d}\Big\| \nu_{n}(\bb) \Big\| + \Big\| \int_{[\tau_0,\tau_j)} \tilde\nu_{n}(\hat\beta(u))dH(u) \Big\| + \Big\| \frac{\tau_0}{n}\sum_{i=1}^n (\ZZ_i - \E\ZZ_i) \Big\| &\leq& \alpha_1
\\
\Big\|\int_{[\tau_0,\tau_j)} \tilde \mu_k(\hat\beta(u)) - \tilde \mu_k(\beta(u)) dH(u) \Big\| &\leq& \alpha_2.
\eean
Then any minimizer of $H_j$ defined in (\ref{intesteq.ava}) is of the form $\Pp_{J}^{-1}(\hat\hh(\tau_j)^t,\vec 0^t)^t$ where $\hat\hh(\tau_j)$ is a minimizer of $H_j(\Pp_{J}^{-1}(\hh^t,\vec 0^t)^t)$ over $\hh\in\R^{|J|}$.
\end{lemma}
\textbf{Proof}
In order to simplify the presentation, assume w.o.l.g. that $J = \{1,...,L\}$, that $\inf_{k\in J} p_k(n,\tau_j) = p_L(n,\tau_j)$ and that $\sup_{k\in J} p_k(n,\tau_j) = p_{L+1}(n,\tau_j)$. Define
\bea \label{esteq.ava2}
\quad\quad \Psi_j(\bb,\xi) &:=& -2\xi^t\Big( \mu(\bb) - \mu(\beta(\tau_{j})) +
\nu_{n}(\bb) - \int_{[\tau_0,\tau_j)} \tilde \nu_{n}(\hat\beta(u))dH(u) \Big)
\\ \nonumber
&& +2\xi^t\int_{[\tau_0,\tau_j)} \tilde\mu(\hat \beta(u)) - \tilde\mu(\beta(u))dH(u) + \frac{1}{n}\sum_{i=1}^n I\{X_i = \ZZ_i^t\bb\}(\delta_i\xi^t\ZZ_i + |\xi^t\ZZ_i|)
\\ \nonumber
&& + 2\lambda_n\sum_{k=1}^d \Big(\xi_k\frac{sgn(\bb_k)}{p_k(n,\tau_j)} + I\{b_k=0\}\frac{|\xi_k|}{p_k(n,\tau_j)}\Big).
\eea
and note that finding all minimizers of the function $H_j((\hh^t,\vec 0^t)^t)$ in (\ref{intesteq.ava}) over $\hh\in\R^L$ is equivalent to finding all points $\hat\bb = (\hat\hh^t,\vec 0^t)^t$ that satisfy
\[
\inf_{\xi = (\zeta^t,\vec 0^t)^t,\zeta\in\R^L} \Psi_j(\hat \bb,\xi) \geq 0.
\]
For a proof of the first part of the lemma, observe that by simple algebraic manipulations and the condition on $\Psi$ we have
\[
0 \leq \Psi_j(\hat \bb,-e_k) = - \Psi_j(\hat \bb,e_k) + \frac{2}{n}\sum_{i=1}^n I\{X_i = \ZZ_i^t\bb\}|e_k^t\ZZ_i| + \frac{4\lambda_n}{p_k(n,\tau_j)}I\{\bb_k=0\}.
\]
This directly yields,
\[
\Psi_j(\hat \bb,e_k) \leq \frac{2}{n}\sum_{i=1}^n I\{X_i = \ZZ_i^t\bb\}|e_k^t\ZZ_i| + \frac{4\lambda_n}{p_k(n,\tau_j)}I\{\bb_k=0\},
\]
and by assumption we have $0 \leq \Psi_j(\hat \bb,e_k)$. From that we obtain for $k=1,...,L$
\bean
&&\Big|\mu_k(\check\beta) - \mu_k(\beta(\tau_j)) + \nu_{n,k}(\check\beta) - \int_{[\tau_0,\tau_j)} \tilde\nu_{n,k}(\hat\beta(u))dH(u) - \int_{[\tau_0,\tau_j)} \tilde \mu_k(\hat\beta(u)) - \tilde \mu_k(\beta(u)) dH(u)
\\
&& \quad - \frac{\tau_0}{n}\sum_{i=1}^n (\ZZ_{i,k} - \E\ZZ_{i,k})\Big|
\\
&=& \frac{1}{2}\Big|\Psi_j(\hat \bb,e_k) - \frac{1}{n}\sum_{i=1}^n I\{X_i = \ZZ_i^t\bb\}(\delta_i\ZZ_{i,k} + |\ZZ_{i,k}|) - \frac{2\lambda_n}{p_k(n,\tau_j)}\Big( sgn(\bb_k)+ I\{b_k=0\}\Big) \Big|
\\
&\leq& \frac{1}{2} \Big(\Big|\Psi_j(\hat \bb,e_k) - \frac{1}{n}\sum_{i=1}^n I\{X_i = \ZZ_i^t\bb\}|\ZZ_{i,k}| - \frac{2\lambda_nI\{\bb_k=0\}}{p_k(n,\tau_j)} \Big| + \frac{2\lambda_nI\{\bb_k\neq 0\}}{p_k(n,\tau_j)} + \frac{C_Z}{n}\Big)
\\
&\leq& \frac{\lambda_n}{p_k(n,\tau_j)} + \frac{C_Z}{n}
\eean
almost surely. A simple application of the triangle inequality completes the proof of the first part of the lemma.\\
For a proof of the second part, assume w.o.l.g. that $ J = \{1,...,L\}$ and that the assumptions made at the beginning of the proof of the first part hold. In particular, under this simplifying assumptions $\Pp_{\tau_j}$ is the identity matrix. Start by noting that
\bean
\Psi_j(\bb,\xi_1+\xi_2) &=& \Psi_j(\bb,\xi_1) + \Psi_j(\bb,\xi_2) - \frac{1}{n}\sum_{i=1}^nI\{X_i = \ZZ_i^t b\}(|\xi_1^t\ZZ_i| + |\xi_2^t\ZZ_i| - |(\xi_1+\xi_2)^t\ZZ_i|)
\\
&& - 2\lambda_n \sum_{k=1}^d \frac{I\{b_k=0\}}{p_k(n,\tau_j)}(|\xi_{1,k}|+|\xi_{2,k}| - |\xi_{1,k}+\xi_{2,k}|).
\eean
In particular, for the special case $\xi_1^0 = (\zeta^t,\vec 0_{d-L}^t)^t, \xi_2^0 = (\vec 0_L^t,\theta^t)^t$ with $\zeta\in\R^L,\theta\in\R^{d-L}$, the last line in the above equation equals zero. Moreover, $|a|+|b|-|a+b|\leq 2|b|$, and thus $|\xi_1^t\ZZ_i| + |\xi_2^t\ZZ_i| - |(\xi_1+\xi_2)^t\ZZ_i| \leq 2|\xi_2^t\ZZ_i|$. Hence, if we can show that
\[
\Psi_j(\check\beta,\xi_1^0) + \Psi_j(\check\beta,\xi_2^0) \geq \frac{2C_Z}{n}\sum_{l=L+1}^{d}|\xi_{2,l}^0|
\]
for any $\xi_1^0,\xi_2^0$ of the form given above, it will follow that $\Psi_j(\check\beta,\xi) \geq 0$ for all $\xi \in \R^d$. By the definition of $\check\beta$ we have $\Psi_j(\check\beta,\xi_1^0) \geq 0$, and thus it remains to verify that $\Psi_j(\check\beta,\xi_2^0) \geq \frac{2C_Z}{n}\sum|\xi_{2,l}^0|$. To this end, observe that the arguments in the first part of the Lemma yield the bound [the last inequality follows under (\ref{bednull.ava})]
\begin{eqnarray*}
\|\bar\mu^{( J)}(\check\beta) - \bar\mu^{( J)}(\beta^{( J)}(\tau_{j}))\| &\leq& \alpha_1 + \alpha_2 + \frac{\lambda_n}{p_k(n,\tau_j)} + \frac{C_Z}{n} + C_2 \sup_{k > L} |\beta_k(\tau_j))|
\\
&\leq& \frac{\eps - \sup_{k > L} |\beta_k(\tau_j))|}{C_1\vee 1}
\end{eqnarray*}
since by assumption (C3), condition (\ref{bednull.ava}) and Remark \ref{rem.mu} we have
\[
\|\bar\mu^{( J)}(\beta(\tau_j)) - \bar\mu^{( J)}(\beta^{( J)}(\tau_j))\| \leq C_2 \sup_{k > L} |\beta_k(\tau_j))|.
\]
Thus $\|\check\beta^{(\chi(\tau_j))} - \beta(\tau_j)\| \leq \eps$ and (\ref{muineq1}) together with the triangle inequality implies that
\[
\|\mu(\check\beta) - \mu(\beta(\tau_{j}))\| \leq C_5\Big(\alpha_1 + \alpha_2 + \frac{\lambda_n}{p_L(n,\tau_j)} + \frac{C_Z}{n}
+ C_2 \sup_{k > L} |\beta_k(\tau_j))| \Big)+ C_2 \sup_{k > L} |\beta_k(\tau_j))|.
\]
By the definition of $\check\beta$ and the assumption on $p_k(n,\tau_j)$ made at the beginning of the proof we have
\[
2\lambda_n\sum_{k=1}^d \Big(\xi^0_{2,k}\frac{sgn(\check\beta_k)}{p_k(n,\tau_j)} + I\{\check\beta_k=0\}\frac{|\xi^0_{2,k}|}{p_k(n,\tau_j)}\Big) = 2\lambda_n\sum_{k=L+1}^d\frac{|\xi^0_{2,k}|}{p_k(n,\tau_j)} \geq \frac{2\lambda_n}{p_{L+1}(n,\tau_j)}\sum_{k=L+1}^{d}|\xi_{2,k}^0|.
\]
Combining all the inequalities derived above, we see from the definition of $\Psi$ that
\bean
\Psi_j(\bb,\xi_2^0) &\geq& \sum_{k=L+1}^{d}|\xi_{2,k}^0|\Big(\frac{2\lambda_n}{p_{L+1}(n,\tau_j)} - 2\alpha_1 - 2\alpha_2 - 2\|\mu(\check\beta) - \mu(\beta(\tau_{j}))\| - \frac{2C_Z}{n}\Big).
\eean
Thus under (\ref{bednull2.ava}) it holds that $\Psi_j(\check\beta,\xi_2^0) \geq \frac{2(C_5+1)C_Z}{n}\sum|\xi_{2,l}^0| \geq \frac{2C_Z}{n}\sum|\xi_{2,l}^0|$ and we have proved that $\check\beta$ is a minimizer of the function $H(\bb)$ in the set $\R^d$. It remains to verify that every minimizer is of this form. We will prove this assertion by contradiction. Assume that there exists a minimizer $\check \bb$ with $\check \bb_k \neq 0$ for some $k>L$. Since the set of minimizers is convex, any convex combination of $\check \bb$ and a minimizer $\check\beta$ with $\check\beta_k=0$ would also be a minimizer. Thus there must exist a minimizer $\tilde\bb$ with k'th component different from zero and all other components arbitrarily close to the components of $\check \beta$. In particular, we can choose $\tilde\bb$ in such a way that $\|\mu(\tilde\bb) -\mu(\check\beta)\| \leq \frac{C_5C_Z}{n}$. Setting $\bb = \tilde\bb, \xi = \pm e_k$ in representation (\ref{esteq.ava2}) we obtain a contradiction, since in this case the sum in the last line will take the values $\pm 2\lambda_n\frac{sgn(\tilde\bb_k)}{p_k(n,\tau_j)}$, and the absolute value of this quantity dominates the rest of $\Psi_j(\tilde\bb,\xi)$ by construction and condition (\ref{bednull2.ava}). Thus a minimizer with $\check \bb_k \neq 0$ for some $k>L$ can not exist and proof is complete. \hfill $\Box$

\newpage

\begin{lemma} \label{lem:consistent}
Under assumptions (C1)-(C4) and (D1) the unpenalized estimators obtained from minimizing (\ref{intesteq.ava_np}) are uniformly consistent in probability, i.e.
\[
\sup_{\tau\in[\tau_L,\tau_U]} \|\hat\beta(\tau) - \beta(\tau)\|= o_P(1).
\]
\end{lemma}

\textbf{Proof}
Define the quantities
\[
R_{n,1} :=  C_M\Big(\sup_{\bb \in \R^d} \|\nu_n(\bb)\| + H(\tau_U) \sup_{\bb \in \R^d}\|\tilde \nu_n(\bb)\| + \Big\|\frac{\tau_0}{n}\sum_{i=1}^n (\ZZ_i-\E\ZZ_i) \Big\|\Big) = o_P(1),
\]
$r_{n,1} := C_5\Big(R_{n,1} + C_6b_n^2 + \frac{2C_Z}{n}\Big)$ and 
\[
\Rr_n := \Big(r_{n,1} +  \frac{C_6b_n}{\tilde C_5}\Big)\sup_n(1+\tilde C_5b_n)^{N_\tau(n)} = o_P(1).
\]
Use similar arguments as in step 1 of the proof of Lemma \ref{lem:dmuphi} [set $\Lambda_{n,1} = 0, \Lambda_{n,1} = +\infty$] to inductively show that on the set $\Omega_n := \Big\{\Rr_n \leq \frac{\eps}{C_1\wedge 1}\Big\}$ whose probability tends to one we have
\begin{enumerate}
\item[(i)] the conditions (\ref{bednull.ava}) and (\ref{bednull2.ava}) of Lemma \ref{lem:bednull.gen} hold with $J = \{1,...,d\}$.
\item[(ii)] we have the following upper bound
\bean
\|\mu(\hat \beta(\tau_j)) - \mu(\beta(\tau_j))\| &\leq& r_{n,1}(1+\tilde C_5 b_n)^j + \frac{C_6b_n}{\tilde C_5}((1+\tilde C_5 b_n)^j-1) =: r_{n,j+1}
\\
&\leq& \Big(r_{n,1} +  \frac{C_6b_n}{\tilde C_5}\Big)\sup_n(1+\tilde C_5b_n)^{N_\tau(n)} = \Rr_n = o_P(1).
\eean
\end{enumerate}
In particular, the results above and an application of Remark \ref{rem.mu} imply that
\[
\sup_{j=1,...,N(\tau)} \|\hat\beta(\tau_j) - \beta(\tau_j)\| = o_P(1).
\]
Since $\hat\beta(\tau)$ is constant between grid points and additionally $\beta(\tau)$ is Lipschitz-continuous, this completes the proof. 

\hfill $\Box$

\newpage

\begin{lemma} \label{lem:dmuphi}
Define the triangular array of random $\R^d$-valued vectors $\phi_n(\tau_j)$ as
\beq \label{phi0}
\phi_n(\tau_0) - \mu(\beta(\tau_{0}))= -\M_{\tau_0}\Big(\frac{1}{n} \sum_{i=1}^n \ZZ_i(I\{X_i \leq \ZZ_i^t\beta(\tau_0)\} - \tau_0)\Big)
\eeq
and for $j= 1,...,N_\tau$
\bea \label{phirek}
\phi_n(\tau_{j}) - \mu(\beta(\tau_{j})) &=& \M_{\tau_j}\Big( -\nu_n(\beta(\tau_j)) + \int_{[\tau_0,\tau_j)} \tilde \nu_n(\beta(u))dH(u)
+ \frac{\tau_0}{n} \sum_{i=1}^n (\ZZ_i - \E\ZZ_i) \quad
\\ \nonumber
&& \quad + \sum_{i=0}^{j-1} \int_{[\tau_i,\tau_{i+1})}\tilde\M_{u} dH(u)\Big(\phi_n(\tau_i) - \mu(\beta(\tau_i))\Big)\Big).
\eea
\begin{enumerate}
\item[(a)] Let assumptions (C1)-(C4) and (D1)-(D3) hold and denote by $\hat \beta$ the unpenalized estimator obtained from minimizing (\ref{intesteq.ava_np}). Then 
\[
\sqrt n \sup_{j} \|\mu(\hat\beta(\tau_j)) - \phi_n(\tau_j)\| = O_P(n^{1/2}b_n + n^{-\gamma/2} + \omega_{c_nn^{-1/2}}(\sqrt{n}\nu_n) + \omega_{c_nn^{-1/2}}(\sqrt{n}\tilde\nu_n))
\]
\item[(b)]
Let assumptions (C1)-(C3), (C4'), (D1'), (D2)-(D3), (P) hold and denote by $\hat \beta$ the penalized estimator obtained from minimizing (\ref{intesteq.ava}) . 
Then $\sqrt n \sup_{j} \|\mu(\hat\beta(\tau_j)) - \phi_n(\tau_j)\| = o_P(1)$ and $P(\sup_{\tau_j}\sup_{k \in \chi(\tau_j)^C} |\hat\beta_k(\tau_j)| = 0) \rightarrow 1$.
\end{enumerate}
\end{lemma}
\textbf{Proof.}
The proof of part (a) is similar to, but simpler then the proof of part (b). For this reason, we will only state the proof of (b) and point out the important differences where necessary. The proof will consist of two major steps. In the first step we define the set
\[
\Omega_n := \Big\{\Rr_n \leq \frac{\eps}{C_1\wedge 1}\Big\} \cap \Big\{(1+C_5)\Rr_n + C_5\Lambda_{n,1} \leq \Lambda_{n,0}\Big\} \cap \Omega_{0,n}
\]
with $\Omega_{0,n}$ denoting some set such that $P(\Omega_{0,n}) \to 1$ and note that $P(\Omega_n) \to 1$, here [the bound will be proved below]
\[
\Rr_n := \Big(r_{n,1} +  \frac{C_6b_n}{\tilde C_5}\Big)\sup_n(1+\tilde C_5b_n)^{N_\tau(n)} = O_P(n^{-1/2})
\]
and $r_{n,1} := C_5\Big(R_{n,1} + C_6b_n^2 + \frac{2C_Z}{n} + \Lambda_{n,1}\Big)$ with
\beq \label{defrn1}
R_{n,1} :=  C_M\Big(\sup_{\bb \in \R^d} \|\nu_n(\bb)\| + H(\tau_U) \sup_{\bb \in \R^d}\|\tilde \nu_n(\bb)\| + \Big\|\frac{\tau_0}{n}\sum_{i=1}^n (\ZZ_i-\E\ZZ_i) \Big\|\Big) = O_P(1/\sqrt n).
\eeq
For a proof of (a), proceed in a similar fashion but with $\chi(\tau)=\{1,...,d\}$ for all $\tau$, setting $\Lambda_{n,1} = 0, \Lambda_{n,0} = \infty$ and replacing $R_{n,1}$ in the definition above by
\[
\tilde R_{n,1} :=  C_M\Big(\sup_{\bb \in \B([\tau_L,\tau_U],\eps)} \|\nu_n(\bb)\| + H(\tau_U) \sup_{\bb \in \B([\tau_L,\tau_U],\eps)}\|\tilde \nu_n(\bb)\| + \Big\|\frac{\tau_0}{n}\sum_{i=1}^n (\ZZ_i-\E\ZZ_i) \Big\|\Big).
\]
Here, uniform consistency of the unpenalized estimator [see Lemma \ref{lem:consistent}] implies that only the supremum over $\bb \in \B([\tau_L,\tau_U],\eps)$ needs to be considered.\\
In what follows, we will inductively show that on the set $\Omega_n$ we have for every $0 \leq j \leq N_\tau(n)$ [recall that $N_\tau(n)$ is the number of grid points]
\begin{enumerate}
\item[(i)] the conditions (\ref{bednull.ava}) and (\ref{bednull2.ava}) of Lemma \ref{lem:bednull.gen} hold [the quantities $\alpha_1,\alpha_2$ will depend on $j$ and be specified in the proof below].
\item[(ii)] $\hat\beta_k(\tau_j)=0$ for $k \in \chi(\tau_j)^C$.
\item[(iii)] we have the following upper bound
\bean
\|\mu^{(\chi(\tau_j))}(\hat \beta(\tau_j)) - \mu^{(\chi(\tau_j))}(\beta(\tau_j))\| &\leq& r_{n,1}(1+\tilde C_5 b_n)^j + \frac{C_6b_n}{\tilde C_5}((1+\tilde C_5 b_n)^j-1) =: r_{n,j+1}
\\
&\leq& \Big(r_{n,1} +  \frac{C_6b_n}{\tilde C_5}\Big)\sup_n(1+\tilde C_5b_n)^{N_\tau(n)} = \Rr_n = O_P(n^{-1/2})
\eean
\end{enumerate}

\noindent
In the second step, we will prove the bounds
\bea \label{vn}
&&  \label{vn1}
\sup_j \|\mu(\hat\beta(\tau_j)) - \phi_n(\tau_j)\| \leq s_{n,1}\sup_n (1+dC_Mb_n)^{N_\tau(n)}
\eea
where $s_{n,1} = o_P(n^{-1/2})$ in case (b) and 
\[
s_{n,1} = O_P(b_n + n^{-(1+\gamma/2)} + \omega_{c_nn^{-1/2}}(\nu_n) + \omega_{c_nn^{-1/2}}(\tilde\nu_n))
\] 
in case (a).
\\
\noindent
\textbf{Step 1: Proof of (i), (ii) and (iii).}\\
First, consider the grid point $\tau_0$. Classical arguments yield the existence of a set $\Omega_{0,n}$ such that $P(\Omega_{0,n})\to 1$ and (ii)-(iii) hold on this set. The details are omitted for the sake of brevity.\\ 
Next, observe that for the grid point $\tau_1$ we have for $k \in \{1,...,d\}$ [apply Remark \ref{rem.mu}]
\[
\Big|\int_{[\tau_0,\tau_{1})} \tilde \mu_k(\hat\beta(u)) - \tilde \mu_k(\beta(u)) dH(u) \Big| \leq r_{n,1} + C_6b_n^2 =: R_{n,2} = O_P(n^{-1/2}).
\]
Defining $\alpha_j := R_{n,j}$ ($j=1,2$) we obtain that conditions (\ref{bednull.ava}) and (\ref{bednull2.ava}) of Lemma \ref{lem:bednull.gen} hold with $j=1$ on the set
\[
\Omega_{1,n} := \Big\{\frac{C_Z}{n} + R_{n,1} + R_{n,2} + \Lambda_{n,1} \leq \frac{\eps}{C_1\wedge 1} \Big\}\cap \Big\{(1+C_5)(\frac{2C_Z}{n} + R_{n,1} + R_{n,2}) + C_5\Lambda_{n,1} \leq \Lambda_{n,0} \Big\}.
\]
Finally, note that by the first part of Lemma \ref{lem:bednull.gen} we have for $k\in\chi(\tau_1)$
\bean
|\mu_k(\hat\beta(\tau_1)) - \mu_k(\beta(\tau_1))| \leq R_{n,1} + R_{n,2} + \frac{2C_Z}{n} + \Lambda_{n,1}
\eean
[the constant $2$ in front of $C_Z$ will play a role later] which implies (iii) on the set $\Omega_{1,n}$.
\\
Now, proceed inductively. Assume that (i)-(iii) have been established for $1,...,j$. For the grid point $\tau_{j+1}$, observe that for $k \in \{1,...,d\}$
\[
\Big|\int_{[\tau_0,\tau_{j+1})} \tilde \mu_k(\hat\beta(u)) - \tilde \mu_k(\beta(u)) dH(u) \Big|
\leq R_{n,2} + b_n\sum_{i=1}^j (\tilde C_5 r_{n,i} + C_6b_n).
\]
Thus, setting $\alpha_1 = R_{n,1}$, $\alpha_2 :=  R_{n,2} + b_n\sum_{i=1}^j (\tilde C_5 r_{n,i} + C_6b_n)$ we obtain that conditions (\ref{bednull.ava}) and (\ref{bednull2.ava}) of Lemma \ref{lem:bednull.gen} hold on the set
\bean
\Omega_{j+1,n} &:=& \Big\{\frac{C_Z}{n} + R_{n,1} + R_{n,2} + b_n\sum_{i=1}^j (C_6b_n + \tilde C_5 r_{n,i}) + \Lambda_{n,1} \leq \frac{\eps}{C_1\wedge 1} \Big\}\cap
\\
&& \quad \quad \cap \Big\{(1+C_5)\Big(\frac{2C_Z}{n} + R_{n,1} + R_{n,2} + b_n\sum_{i=1}^j (C_6b_n + \tilde C_5 r_{n,i})\Big) + C_5\Lambda_{n,1} \leq \Lambda_{n,0} \Big\}.
\eean
This yields (i) and (ii) for $\tau_{j+1}$ on the set $\Omega_{j+1,n}$. Finally, note that by the first part of Lemma \ref{lem:bednull.gen} we have for $k \in \chi(\tau_j)$
\bean
|\mu_k(\hat\beta(\tau_{j+1})) - \mu_k(\beta(\tau_{j+1}))| \leq r_{n,1} + b_n\sum_{i=1}^j (\tilde C_5 r_{n,i} + C_6b_n).
\eean
Inserting the definition of $r_{n,k}$ for $k=2,...,j$, some algebra yields
\[
r_{n,1} + b_n\sum_{i=1}^j (\tilde C_5 r_{n,i} + C_6b_n)
= r_{n,1}(1+\tilde C_5b_n)^j + C_6b_n\frac{(1+\tilde C_5b_n)^j - 1}{\tilde C_5} = r_{n,j+1},
\]
which completes the proof of (iii) for $\tau_{j+1}$. This shows $\Omega_n\subset\cap_{j}\Omega_{j,n} $ and completes the first step.
\\
\\
\textbf{Step 2:}\\
First of all, note that (iii) from the first step in combination with Remark \ref{rem.mu} shows that 
\bea \label{eq.rate1}
\sup_j \|\hat\beta(\tau_j) - \beta(\tau_j)\| = O_P(n^{-1/2}).
\eea
In order to establish (\ref{vn1}), note that on the set $\Omega_n$ Lemma \ref{lem:bednull.gen} in combination with Remark \ref{rem.mu} yields
\bean
&& \|\mu(\hat\beta(\tau_{j})) - \phi_n(\tau_{j})\| = \|\phi_n(\tau_{j}) - \mu(\beta(\tau_{j})) - (\mu(\hat\beta(\tau_{j})) - \mu(\beta(\tau_{j})))\|
\\
&\leq& \Big\| \M_{\tau_j}\Big( -\nu_n(\beta(\tau_j)) + \int_{[\tau_0,\tau_j)} \tilde \nu_n(\beta(u))dH(u) + \frac{\tau_0}{n} \sum_{i=1}^n (\ZZ_i - \E\ZZ_i)\Big)
\\
&& \quad + \sum_{i=0}^{j-1} \int_{[\tau_i,\tau_{i+1})}\M_{\tau_j}\tilde\M_{u} dH(u)\Big(\phi_n(\tau_i) - \mu(\beta(\tau_i))\Big)
\\
&& \quad + \M_{\tau_j}\Big(\nu_{n}(\hat\beta(\tau_{j})) - \int_{[\tau_0,\tau_j)} \tilde\nu_{n}(\hat\beta(u))dH(u) -\frac{\tau_0}{n} \sum_{i=1}^n (\ZZ_i - \E\ZZ_i)
\\
&& \quad - \int_{[\tau_0,\tau_j)} \tilde \mu(\hat\beta(u)) - \tilde\mu(\beta(u)) dH(u) \Big)\Big\|
\\
&& + \Big\|\mu(\hat\beta(\tau_{j})) - \mu(\beta(\tau_{j})) - \M_{\tau_j}(\bar\mu^{(\chi(\tau_j))}(\hat\beta(\tau_{j})) - \bar\mu^{(\chi(\tau_j))}(\beta(\tau_j)))\Big\|
\\
&& + \Lambda_{n,1} + \frac{C_Z}{n}.
\eean
Now for $n$ large enough and $c_n \to\infty$ we have by (\ref{eq.rate1})
\[
\Big\|- \nu_n(\beta(\tau_j)) + \int_{[\tau_0,\tau_j)} \tilde\nu_n(\beta(u))dH(u) + \Big(\nu_{n}(\hat\beta(\tau_{j})) - \int_{[\tau_0,\tau_j)} \tilde\nu_{n}(\hat\beta(u))dH(u) \Big)\Big\| \leq V_n
\]
where
$
V_n := \omega_{c_nn^{-1/2}}(\nu_n) + H(\tau_U)\omega_{c_nn^{-1/2}}(\tilde\nu_n)
$
and moreover [here, $\tilde\D(\tau)$ is defined in (\ref{def.d})] 
\bean
&&\Big\| \int_{[\tau_j,\tau_{j+1})}\tilde \mu(\hat\beta(u)) - \tilde \mu(\beta(u)) dH(u) - \int_{[\tau_j,\tau_{j+1})}\tilde\M_u dH(u)(\mu(\hat\beta(\tau_j)) - \mu(\beta(\tau_j))) \Big\|
\\
&\leq& (\tilde\D(\Rr_n) + db_nC_MC_7)(H(\tau_{j+1})-H(\tau_{j})).
\eean
In particular, this implies
\bean
&&\Big\|\int_{[\tau_0,\tau_j)}\tilde \mu(\hat\beta(u)) - \tilde \mu(\beta(u)) dH(u) - \sum_{i=0}^{j-1} \int_{[\tau_i,\tau_{i+1})}\tilde\M_udH(u)\Big(\phi_n(\beta(\tau_i)) - \mu(\beta(\tau_i))\Big) \Big\|
\\
&\leq& H(\tau_U)(\tilde\D(\Rr_n) + db_nC_MC_7) + d b_nC_M \sum_{i=0}^{j-1}\|\mu(\hat\beta(\tau_i)) - \phi_n(\tau_i)\|.
\eean
Summarizing, we have obtained that for $j\geq 0$ on the set $\Omega_n$
\bean
&& \|\mu(\hat\beta(\tau_{j})) - \phi_n(\tau_{j})\|
\\
&\leq& \Lambda_{n,1} + \frac{C_Z}{n} + V_n + \D(\Rr_n) + H(\tau_U)(\tilde\D(\Rr_n) + db_nC_MC_7) + d b_nC_M \sum_{i=0}^{j-1}\|\mu(\hat\beta(\tau_i)) - \phi_n(\tau_i)\|.
\eean
Defining
\bean
s_{n,1} &:=& \Lambda_{n,1} + \frac{C_Z}{n} + V_n + \D(\Rr_n) + H(\tau_U)(\tilde\D(\Rr_n) + db_nC_MC_7)
\\
s_{n,j+1} &:=& s_{n,1} + db_n C_M \sum_{i=0}^j s_{n,i}
\eean
we obtain $\|\mu(\hat\beta(\tau_{j+1})) - \phi_n(\tau_{j+1})\| \leq s_{n,j+1}$. Moreover, induction yields
\[
s_{n,j+1} = (1+dC_Mb_n)^{j+1}s_{n,1} \leq s_{n,1}\sup_n (1+dC_Mb_n)^{N_\tau(n)}. 
\]
This completes the proof.
\hfill $\Box$

\newpage

\begin{lemma} \label{lem:phi}
Under the assumptions of Lemma \ref{lem:dmuphi} we have for $j=0,...,N_\tau(n)$
\bean
&&\phi_n(\tau_j) - \mu(\beta(\tau_j))
\\
&=& \M_{\tau_j}\Big( w_n(\tau_j) + \int_{[\tau_0,\tau_j)}\Big(\prodi_{(u,\tau_j]}\Big(I_{d} + (\M_v \tilde\M_v)^tdH(v) \Big)\Big)^t \tilde\M_u\M_u w_n(u)dH(u)\Big)
+ R_n(\tau_j)
\eean
uniformly in $j$ where for some finite constant $C$ we have $\sup_j \|R_n(\tau_j)\| = O_P(r_n)$ with
\bean
r_n &:=& C\Big(\Big(b_n + \sup_{|u-v|\leq a_n, \theta_k \notin [u,v] \forall k}\|\M_u - \M_v\| + \|\tilde\M_u - \tilde\M_v\|\Big)\sup_u \|w_n(u)\| 
\\
&& + \sup_{|u-v| \leq a_n} \|w_n(u)-w_n(v)\|\Big),
\eean
$I_{d}$ denotes the $d\times d$ identity matrix, $\prodi$ denotes the product-integral [see \cite{gilljoha1990}] and we defined
\[
w_n(\tau) := \frac{\tau_0}{n}\sum_{i=1}^n (\ZZ_i - \E\ZZ_i)- \nu_n(\tau) + \int_{[\tau_0,\tau)}\tilde \nu_n(u)dH(u).
\]
\end{lemma}
\textbf{Proof.}
Throughout this proof, denote by $C$ some generic constant whose value might differ from line to line. Start by noting that the solution of the iterative equation (\ref{phirek}) is given by
\bean
\phi_n(\tau_{j+1}) - \mu(\beta(\tau_{j+1})) &=& \M_{\tau_{j+1}}\sum_{l=0}^j \Big(\prod_{i=l+1}^{j}\Big(I_d + \int_{[\tau_i,\tau_{i+1})} (\M_{\tau_i}\tilde\M_u)^t dH(u) \Big)\Big)^t(w_n(\tau_{l+1})-w_n(\tau_l))
\\
&& + \M_{\tau_{j+1}}\Big(\prod_{i=0}^{j}\Big(I_d + \int_{[\tau_i,\tau_{i+1})} (\M_{\tau_i}\tilde\M_u)^t dH(u)\Big)\Big)^tw_n(\tau_0),
\eean
this assertion can be proved by induction [here, a product $\prod_{i=a}^b C_i$ with $a>b$ is defined as the unit matrix of suitable dimension]. Next, observe that summation-by-parts, that is
\[
\sum_{k=m}^n f_k(g_{k+1}-g_k) = f_{n+1}g_{n+1} - f_mg_m - \sum_{k=m}^n (f_{k+1}-f_k)g_{k+1}
\]
yields 
\bean
&& \sum_{l=0}^j \Big(\prod_{i=l+1}^{j}\Big(I_d + \int_{[\tau_i,\tau_{i+1})}(\M_{\tau_i}\tilde\M_u)^tdH(u) \Big)\Big)^t(w_n(\tau_{l+1})-w_n(\tau_l))
\\
&& \quad\quad\quad\quad + \Big(\prod_{i=0}^{j}\Big(I_d + \int_{[\tau_i,\tau_{i+1})} (\M_{\tau_i}\tilde\M_u)^t dH(u)\Big)\Big)^tw_n(\tau_0)
\\
&=& I_d w_n(\tau_{j+1})  - \sum_{l=0}^j \Big[\prod_{i=l+2}^{j}\Big(I_d + \int_{[\tau_i,\tau_{i+1})}(\M_{\tau_i}\tilde\M_u)^tdH(u) \Big)
\\
&& \quad\quad\quad\quad - \prod_{i=l+1}^{j}\Big(I_d +\int_{[\tau_i,\tau_{i+1})}(\M_{\tau_i}\tilde\M_u)^tdH(u) \Big)\Big]^t w_n(\tau_{l+1})
\\
&=& w_n(\tau_{j+1}) + \sum_{l=0}^{j-1} \Big( \prod_{i=l+2}^{j}\Big(I_d + \int_{[\tau_i,\tau_{i+1})}(\M_{\tau_i}\tilde\M_u)^tdH(u) \Big)  \Big)^t \int_{[\tau_{l+1},\tau_{l+2})}\tilde\M_u\M_{\tau_{l+1}} dH(u)w_n(\tau_{l+1}).
\eean
At the end of the proof, we will show that
\beq \label{estprodi}
\quad \ \ \sup_{k,j,j<k} \Big\|\prod_{i=j}^{k-1}\Big(I_d + \int_{[\tau_i,\tau_{i+1})}(\M_{\tau_i}\tilde\M_u)^tdH(u) \Big) -  \prodi_{(\tau_j,\tau_k]}\Big(I_d + (\M_u\tilde\M_u)^tdH(u) \Big) \Big\| \leq Cd_n
\eeq
where $d_n := b_n + \sup_{|u-v|\leq a_n, \theta_k \notin [u,v] \forall k}\Big(\|\M_u - \M_v\| + \|\tilde\M_u - \tilde\M_v\|\Big)$.
Moreover, we note that
\bean
&& \sup_{k,j,j<k} \sup_{v \in (\tau_j,\tau_{j+1}]} \Big\|\prodi_{(\tau_j,\tau_{k}]}\Big(I_d + (\M_u\tilde\M_u)^tdH(u) \Big)-\prodi_{(v,\tau_{k}]}\Big(I_d + (\M_u\tilde\M_u)^tdH(u) \Big) \Big\| \leq Cb_n
\eean
since $\Big\|\prodi_{(a,b]}\Big(I_d + (\M_u\tilde\M_u)^tdH(u) \Big)\Big\| \leq \exp(dC_M(H(b)-H(a)))$ by inequality (37) from \cite{gilljoha1990} and $\Big\|\prodi_{(v,\tau_{j+1}]}\Big(I_d + (\M_u\tilde\M_u)^tdH(u) \Big) - I_d\Big\| \leq dC_M(H(\tau_{j+1})-H(v))\exp(dC_M(H(\tau_{j+1})-H(v)))$ by inequality (38) from the same reference. This yields
\bean
&& \sup_j\Big\| \sum_{l=0}^{j-1} \Big( \prod_{i=l+2}^{j}\Big(I_d + \int_{[\tau_i,\tau_{i+1})}(\M_u\tilde\M_u)^tdH(u)\Big)\Big)^t \int_{[\tau_{l+1},\tau_{l+2})}(\M_u\tilde\M_u)^tdH(u)w_n(\tau_{l+1})
\\
&& \quad\quad\quad\quad - \int_{[\tau_0,\tau_{j+1})}\Big(\prodi_{(v,\tau_{j+1}]}\Big(I_d + (\M_u\tilde\M_u)^tdH(u) \Big)\Big)^t \tilde\M_v \M_v w_n(v)dH(v) \Big\| \leq C r_n,
\eean
since $\int_{[\tau_0,\tau_{1})}\Big(\prodi_{(v,\tau_{1}]}\Big(I_d + (\M_u\tilde\M_u)^tdH(u) \Big)\Big)^t \tilde\M_v\M_v w_n(v)dH(v) \leq C b_n\sup_u \|w_n(u)\|$. Thus it remains to establish (\ref{estprodi}). To this end, we note that
\bean
&& \prod_{i=j}^{k-1}\Big(I_d + \int_{[\tau_i,\tau_{i+1})}(\M_{\tau_i}\tilde\M_u)^tdH(u) \Big) -  \prod_{i=j}^{k-1}\Big(I_d + \int_{[\tau_i,\tau_{i+1})}(\M_{u}\tilde\M_u)^tdH(u) \Big)
\\
&=& \sum_{l=j}^{k-1} \Big(\prod_{i=j}^{l-1}\Big(I_d + \int_{[\tau_i,\tau_{i+1})}(\M_{\tau_i}\tilde\M_u)^tdH(u) \Big)\Big)
\\
&&\quad\quad\times
\Big(\int_{[\tau_l,\tau_{l+1})}(\M_{\tau_l}\tilde\M_u)^tdH(u) - \int_{[\tau_l,\tau_{l+1})}(\M_{u}\tilde\M_u)^tdH(u) \Big)\Big)
\\
&&\quad\quad\times \prod_{i=l+1}^{k-1}\Big(I_d + \int_{[\tau_i,\tau_{i+1})}(\M_{u}\tilde\M_u)^tdH(u) \Big).
\eean
Next, observe that
\[
\sup_{l:\theta_k \notin [\tau_l,\tau_{l+1})\forall k}\Big\| \int_{[\tau_l,\tau_{l+1})}(\M_{\tau_l}\tilde\M_u)^tdH(u) - \int_{[\tau_l,\tau_{l+1})}(\M_u\tilde\M_u)^tdH(u) \Big\| \leq  C b_n  \sup_{|u-v|\leq a_n, \theta_k \notin [u,v] \forall k}\|\M_u - \M_v\| 
\]
and
\[
\sup_{l:\exists k: \theta_k \in [\tau_l,\tau_{l+1})}\Big\| \int_{[\tau_l,\tau_{l+1})}(\M_{\tau_l}\tilde\M_u)^tdH(u) - \int_{[\tau_l,\tau_{l+1})}(\M_u\tilde\M_u)^tdH(u) \Big\| \leq Cb_n.
\]
Finally, note that $k-1-j \leq N_\tau(n)$, $b_nN_\tau(n) = O(1)$ and that $\Big\|\int_{[\tau_l,\tau_{l+1})}\M_{\tau_l}\tilde\M_udH(u)\Big\| \leq Cb_n$, $\Big\|\int_{[\tau_l,\tau_{l+1})}\M_{\tau_u}\tilde\M_udH(u)\Big\| \leq Cb_n$ uniformly in $l$, which yields
\[
\sup_{k,j}\Big\|
\prod_{i=j}^{k-1}\Big(I_d + \int_{[\tau_i,\tau_{i+1})}(\M_{\tau_i}\tilde\M_u)^tdH(u) \Big) -  \prod_{i=j}^{k-1}\Big(I_d + \int_{[\tau_i,\tau_{i+1})}(\M_{u}\tilde\M_u)^tdH(u) \Big)
\Big\| 
\leq Cd_n
\]
since there are only finitely many different $\theta_k$. Finally, the bound
\[
\sup_{k,j,j<k} \Big\|\prod_{i=j}^{k-1}\Big(I_d + \int_{[\tau_i,\tau_{i+1})}(\M_{u}\tilde\M_u)^tdH(u) \Big) -  \prodi_{(\tau_j,\tau_k]}\Big(I_d + (\M_u\tilde\M_u)^tdH(u) \Big) \Big\| \leq Cd_n
\]
can be established by using equations (37), (39) in \cite{gilljoha1990} and the representation
\bean
&& \prod_{i=j}^{k-1}\Big(I_d + \int_{[\tau_i,\tau_{i+1})}(\M_{u}\tilde\M_u)^tdH(u) \Big) -  \prod_{i=j}^{k-1} \prodi_{(\tau_i,\tau_{i+1}]}\Big(I_d + (\M_u\tilde\M_u)^tdH(u) \Big)
\\
&=& \sum_{l=j}^{k-1} \Big(\prod_{i=j}^{l-1}\Big(I_d + \int_{[\tau_i,\tau_{i+1})}(\M_{u}\tilde\M_u)^tdH(u) \Big)\Big)\times
\\
&&\quad\quad\times
\Big(I_d + \int_{[\tau_l,\tau_{l+1})}(\M_{u}\tilde\M_u)^tdH(u) - \prodi_{(\tau_{l},\tau_{l+1}]}\Big(I_d + (\M_u\tilde\M_u)^tdH(u) \Big)\Big)\times
\\
&&\quad\quad\times\prodi_{(\tau_{l+1},\tau_k]}\Big(I_d + (\M_u\tilde\M_u)^tdH(u) \Big).
\eean
This completes the proof.
\hfill $\Box$

\newpage

\subsection{Proof of Theorem \ref{bahadur} and Theorem \ref{satz1}} \label{beweissatz1}
The convergence $P(\sup_{\tau\in [\tau_0,\tau_U]} \sup_{k \in \chi(\tau)^C}\|\hat\beta_k(\tau)| = 0) \rightarrow 1$ is a direct consequence of the results in Lemma \ref{lem:dmuphi}. \\
Next, observe that $\sup_j \sup_{u \in (\tau_j,\tau_{j+1}]}\|\mu(\beta(u)) - \mu(\beta(\tau_{j+1}))\| = O(b_n)$ and similarly 
\[
\sup_j \sup_{u \in (\tau_j,\tau_{j+1}]}\|\psi_n(u) - \psi_n(\tau_{j+1})\| = O(b_n \sup_\tau \|w_n(\tau)\| + \omega_{a_n}(w_n)) \quad a.s. 
\] 
where we defined
\[
\psi_n(\tau) := w_n(\tau) - \int_{[\tau_0,\tau)} \Big(\prodi_{(u,\tau_k]}\Big(I_{d} + (\M_v\tilde\M_v)^tdH(v) \Big)\Big)^t \tilde\M_u\M_u w_n(u)dH(u).
\]
Together with the results in Lemma \ref{lem:dmuphi}  and \ref{lem:phi}, this yields the representation
\bean
&&\mu(\hat \beta(s)) - \mu(\beta(s))
\\
&=& \M_s \Big(w_n(s) + \int_{[\tau_0,s)}\Big(\prodi_{(u,s]}\Big(I_{d} + (\M_v\tilde\M_v)^tdH(v) \Big)\Big)^t \tilde\M_u \M_u w_n(u)dH(u)\Big)
+ R_n(s)
\eean
uniformly in $s \in [\tau_0,\tau_U]$ where  
\[
\sup_{\tau \in [\tau_L,\tau_U]} \sqrt{n}\|R_n(\tau)\| = O_P(n^{1/2}b_n + n^{-\gamma/2} + \omega_{c_nn^{-1/2}}(\sqrt{n}\nu_n) + \omega_{c_nn^{-1/2}}(\sqrt{n}\tilde\nu_n))
\]
under the assumptions of Theorem \ref{bahadur} and $\sup_{\tau \in [\tau_L,\tau_U]} \sqrt{n}\|R_n(\tau)\| = o_P(1)$ under the assumptions of Theorem \ref{satz1}. Thus we have obtained representation (\ref{betarep}), and a Taylor expansion combined with some simple algebra yields (\ref{betarep_bah}).\\
The weak convergence statements in both Theorems follow by the continuous mapping theorem [note that by assumption (A3) and equation (37) from \cite{gilljoha1990}, the components of the matrix $\Big(\prodi_{(u,\tau]}\Big(I_{d} + (\M_v\tilde\M_v)^tdH(v) \Big)\Big)^t \tilde\M_u \M_u $ are uniformly bounded], and thus the proof is complete.
\hfill $\Box$

\newpage

\newpage
\subsection{Proof of Theorem  \ref{lem:unifrates} } \label{bewsatz2}

The following result can be proved by similar arguments as Lemma \ref{lem:bednull.gen}.

\begin{lemma}\label{lem:bednull}
Let conditions (C1)-(C3) and (C4*) hold. Assume that $K \subseteq \xi(\tau_j)$ satisfies the following conditions
\bean
\sup_{\bb\in\R^d}\Big\| \nu_{n}(\bb) \Big\| + \Big\| \int_{[\tau_0,\tau_j)} \tilde\nu_{n}(\hat\beta(u))dH(u) \Big\| + \Big\| \frac{\tau_0}{n}\sum_{i=1}^n (\ZZ_i - \E\ZZ_i) \Big\| &\leq& \alpha_1
\\
\Big\|\int_{[\tau_0,\tau_j)} \tilde \mu(\hat\beta(u)) - \tilde \mu(\beta(u)) dH(u) \Big\| &\leq& \alpha_2
\eean
and
\bea
\label{bednull}
&& \quad \alpha_1 + \alpha_2 + \sup_{k \in K}\frac{\lambda_n}{p_k(n,\tau_j)} + \frac{C_Z}{n}
+ C_2 \sup_{k \in K^C} |\beta_k(\tau_j))| \leq \frac{\eps - \sup_{k \in K^C} |\beta_k(\tau_j)|}{C_1}
\\
\label{bednull2}
&& (C_5+1)\Big(\alpha_1 + \alpha_2  + \frac{2C_Z}{n} +
\sup_{k \in K^C} |\beta_k(\tau_j))|\Big) + C_5\sup_{k \in K}\frac{\lambda_n}{p_k(n,\tau_j)}  \leq \inf_{k \in K^C}\frac{\lambda_n}{p_{k}(n,\tau_j)}.\quad
\eea
Then any minimizer of $H_j$ defined in (\ref{intesteq.ava}) is of the form $\Pp_{K}^{-1}(\hat\hh(\tau_j)^t,\vec 0^t)^t$ where $\hat\hh(\tau_j)$ is a minimizer of $H_j(\Pp_{K}^{-1}(\hh^t,\vec 0^t)^t)$ over $\hh\in\R^{|K|}$. Moreover, it holds that
\[
\|\mu^{(K)}(\hat\beta(\tau_j)) - \mu^{(K)}\beta(\tau_j)\| \leq  \frac{C_Z}{n} + \sup_{k \in K}\frac{\lambda_n}{p_k(n,\tau_j)} + C_2\sup_{k \in K^C} |\beta_k(\tau_j)| + \alpha_1 + \alpha_2.
\]
\end{lemma}

For the proof of Theorem \ref{lem:unifrates}, we will consider points $\tau_j$ such that $\tau_j \in \bigcap_k (B_k\cup V_k)$ and $\tau_j \in P \cup S$ separately. Note that for sufficiently large $n$, the set $P \cup S$ is a union of finitely many disjoint intervals. Without loss of generality, assume that $[\tau_0,\tau_{N_1}] \subset \bigcap_k (B_k\cup V_k)$ and $[\tau_{N_1+1},\tau_{N_2}] \subset P \cup S$, $[\tau_{N_2+1},\tau_{N_3}] \subset \bigcap_k (B_k\cup V_k)$ and so on [of course, $N_1,N_2,...$ depend on $n$, but we do not reflect this fact in the notation].\\
Introduce the 'oracle' penalty $p_k^O(n,\tau_j) := \infty I\{\beta_k(\tau_j) = 0\}$ and define $\hat\beta^O(\tau_j)$ as the solution of the minimization in (\ref{intesteq.ava}) based on this penalty. The basic idea for proving process convergence is to show, that the 'estimator' $\hat\beta^O(\tau_j)$ and $\hat\beta(\tau_j)$ have the same first-order asymptotic expansion uniformly on $\tau_j \in P \cup S$. More precisely, we will show that
\beq \label{eq1}
\sup_{\tau_j \in P \cup S} \|\mu(\hat\beta(\tau)) - \mu(\hat\beta^O(\tau)\| = o_P(n^{-1/2}).
\eeq
Note that by the arguments in the proof of Theorem \ref{satz1} this directly implies the weak convergence in (\ref{alawk}).\\
In order to study the uniform rate of convergence of $\hat\beta(\tau_j)$ on $\bigcap_k (B_k\cup V_k)$, we need to introduce some additional notation
Consider the non-overlapping sets
\[
A_{j,n} := \{t: n^{-1/4}\kappa_n^{-1/2}c_n^{-j/5d} \geq t > n^{-1/4}\kappa_n^{-1/2}c_n^{-(j+1)/5d}\}, \quad j=1,...,5d-1.
\]
Observe that for any $\tau$, the components of $\beta(\tau)$ are contained in at most $d$ of those sets and thus  for any $\tau$ there exist three consecutive sets containing no component of $\beta(\tau)$. Moreover, the diameter of each $A_{j,n}$ is by construction of larger order then $n^{-1/2}$. Thus there exists a function $j(\tau)$ such that the probability of the set
\[
\Omega_A := \{ |\tilde \beta_k(\tau)| \notin A_{j(\tau),n}, \ k=1,...,d,\ \tau \in [\tau_L,\tau_U]\}
\]
tends to one. We will use Lemma \ref{lem:bednull} to show that in each step, coefficients with absolute value below $n^{-1/4}\kappa_n^{-1/2}c_n^{-(j(\tau_k)+1)/5d}$ will be set to zero with probability tending to one.\\
Define the quantities
\bea \label{kappaasy1}
M_{n,j} &:=& \sup_{\tau \in B_j} \frac{\lambda_n}{|\tilde \beta_{j}(\tau)|} = o_P(1/\sqrt n)
\\
\label{kappaasy2}
L_{n,j} &:=& \sqrt n\inf_{\tau \in S_{j} \cup V_j} \frac{\lambda_n}{|\tilde \beta_j(\tau)|} \stackrel{P}{\longrightarrow} \infty.
\\
\label{kappaasy3}
W_{n,j} &:=& \sup_{\tau \in P_j \cup B_j} \frac{\lambda_n}{|\tilde \beta_j(\tau)|} = O_P\Big(\frac{c_n}{n^{1/4}\kappa_n^{1/2}} \Big)
\eea
and $M_n := \sup_j M_{n,j}, L_n := \inf_j L_{n,j}, W_n := \sup_j W_{n,j}$.\\
Now begin by considering $\tau_j \in [\tau_{0},\tau_{N_1}]$. A careful inspection of the proofs of Lemma \ref{lem:dmuphi}, Lemma \ref{lem:phi} and Theorem \ref{satz1} show that the arguments and expansions derived there continue to hold and in particular that
\[
\sup_{\tau_j \in [\tau_0,\tau_{N_1}]} \|\mu(\hat\beta(\tau)) - \mu(\hat\beta^O(\tau)\| = o_P(n^{-1/2})
\]
and
\[
R_{n,2} := \Big\| \int_{[\tau_0,\tau_{N_1})} \tilde\mu(\hat \beta(u)) - \tilde\mu(\beta(u))dH(u) \Big\| = O_P(1/\sqrt n).
\]
Next, consider $\tau_j \in [\tau_{N_1+1},\tau_{N_2}]$.
Define the quantities
\bean
U_n &:=& \inf_{\tau \in [\tau_L,\tau_U]}\inf_{j\in\tilde P(\tau), k\in \tilde S(\tau)}\Big\{\frac{\lambda_n}{|\tilde \beta_k(\tau)|} - \frac{C_5\lambda_n}{|\tilde \beta_j(\tau)|} - \frac{C_2(1+C_5)}{n^{1/4}\kappa_n^{1/2}c_n^{(j(\tau)+1)/5d}}\Big\},
\\
s_{n,0} &:=& C_1\Big(\frac{2C_Z}{n}+\frac{C_2}{n^{1/4}\kappa_n^{1/2}c_n^{1/5d}}+W_n+R_{n,1}+R_{n,2}\Big) + \frac{1}{n^{1/4}\kappa_n^{1/2}c_n^{1/5d}},
\eean
where
\bean
\tilde P(\tau)&:=&\Big\{j\in\{1,\dots,d\}:\ |\beta_j(\tau)| \geq \frac{c_n^{-j(\tau)/5d}}{n^{1/4}\kappa_n^{1/2}}\Big\},\\
\tilde S(\tau)&:=&\Big\{j\in\{1,\dots,d\}:\ |\beta_j(\tau)| \leq \frac{c_n^{-(j(\tau)+1)/5d}}{n^{1/4}\kappa_n^{1/2}}\Big\}.
\eean
Note that by the assumptions on $\kappa_n,c_n$ we have that $U_n$ is at least of the order $\lambda_nn^{1/4}\kappa_n^{1/2}c_n^{1/5d}$ which is of larger order then $n^{-1/2}$. In particular, this implies that the probability of the set
\[
\bar \Omega_{2,n} := \Big\{ (1+C_5)\Big(\frac{2C_Z}{n} + R_{n,1} + R_{n,2} + (N_2-N_1)b_nC_3C_{L,1}(s_{n,0} + b_nC_4)\Big)
\leq U_n\Big\}
\]
where $C_{L,1} := \sup_n (1+C_1C_3b_n)^{N_\tau(n)}<\infty$, tends to one by assumption (B3) since $(N_2-N_1)b_n = O(c_n/\kappa_n)^\gamma$. In the following, we will show that on the set
\[
\Omega_{3,n} := \bar \Omega_{2,n} \cap \Omega_A \cap \Big\{\frac{s_{n,0}}{C_1} + (N_2-N_1)b_nC_3C_{L,1}(s_{n,0}+ C_4b_n) \leq \frac{\eps - \sup_{k \in K^C} |\beta_k(\tau_j)|}{C_1} \Big\}
\]
it holds that for $l=0,...,N_2-N_1$
\bean
\Big \|\int_{[\tau_{N_1},\tau_{N_1+l})}\tilde\mu(\hat\beta(u))-\tilde\mu(\beta(u))dH(u) \Big\|
&\leq&  \sum_{j=0}^{l-1} C_3b_n(s_{n,j}+C_4b_n)
\leq lb_nC_3C_{L,1}(s_{n,0}+b_nC_4)
\\
\|\hat\beta(\tau_{N_1+l}) - \beta(\tau_{N_1+l})\| &\leq& s_{n,l},
\eean
where $s_{n,l}$ satisfies the relation
\bean
s_{n,l+1} &=& s_{n,0} + C_1C_3b_n\sum_{i=0}^{l} (s_{n,i} + C_4b_n) = (1+C_1C_3b_n)^{l+1}s_{n,0} + b_n^2C_1C_3C_4\sum_{j=0}^l (1+b_nC_1C_3)^j
\\
&\leq& C_{L,1} s_{n,0} + C_{L,1} b_n C_4-C_4b_n.
\eean
Note that the assertion for $\tilde \mu$ inductively follows from the assertions for $\beta$ and $s_{n,l}$. To establish those assertions, start by considering the case $l=0$. Let $|\beta_j(\tau_{N_1})| \geq n^{-1/4} \kappa_n^{-1/2} c_n^{-(j(\tau_{N_1}))/5d}$ if and only if $j \in K_0$. By construction, conditions (\ref{bednull}) and (\ref{bednull2}) hold on the set $\Omega_{3,n}$ with $K = K_0$, $\alpha_1 = R_{n,1}$ [with $R_n,1$ defined in equation (\ref{defrn1})], $\alpha_2 = R_{n,2}$. Thus on $\Omega_{3,n}$ we have $\hat \beta_k(\tau_{N_1}) = 0$ for $k \in K_0^C$ and by Lemma \ref{lem:bednull} it holds that $\|\hat \beta(\tau_{N_1}) - \beta(\tau_{N_1}) \| \leq  s_{n,0}$. The rest of the assertion follows by iterating the above argument with $\alpha_1 = R_{n,1}$, $\alpha_2 = R_{n,2} + \sum_{j=0}^{l-1} C_3b_n(s_{n,j}+C_4b_n)$ in the $l$'th step.\\
This yields the assertions (\ref{alazero}) and (\ref{alarate}) on the set $[\tau_L,\tau_{N_2}]$. Note by the computations above
\bean
\Big \|\int_{[\tau_{N_1},\tau_{N_2})}\tilde\mu(\hat\beta(u))-\tilde\mu(\beta(u))dH(u) \Big\|
&\leq& (N_2-N_1)b_nC_3C_{L,1}(s_{n,0}+C_4b_n)
\\
&=& O\Big(\frac{c_n}{\kappa_n}\Big)^\gamma O_P\Big(\frac{c_n}{\kappa_n^{1/2}n^{1/4}} \Big) = o_P(1/\sqrt n).
\eean
In particular, this implies that
\[
\Big \|\int_{[\tau_{N_1},\tau_{N_2})}\tilde\mu(\hat\beta(u))-\tilde\mu(\beta(u))dH(u) - \int_{[\tau_{N_1},\tau_{N_2})}\tilde\mu(\hat\beta^ O(u))-\tilde\mu(\beta(u))dH(u)\Big\| = o_P(n^{-1/2}).
\]
Thus we obtain
\[
\sup_{\tau\in[\tau_{N_2+1},\tau_{N_3}]} \|\mu(\hat\beta^O(\tau)) - \mu(\hat\beta(\tau))\| = o_P(n^{-1/2})
\]
by an iterative application of Lemma \ref{lem:bednull}, the arguments are similar to the ones used in the proofs of Lemma \ref{lem:dmuphi}, Lemma \ref{lem:phi} and Theorem \ref{satz1}. Finally, since the set $P\cup S$ is by assumption a finite union of intervals, we can repeat the arguments above to extend the proof to the whole interval $[\tau_L,\tau_U]$. \hfill $\Box$

\bigskip

 \bibliographystyle{apalike}

\bibliography{StansenCensQuantBibalt}

\end{document}